\documentclass[a4paper, 11pt]{article}
\usepackage{srcltx,graphicx}
\usepackage{amsmath, amssymb, amsthm}
\usepackage{color, xcolor}
\usepackage{array}
\usepackage{subfig}
\usepackage{lscape}
\usepackage{algorithm}      
\usepackage[noend]{algpseudocode}
\usepackage{multirow}
\usepackage{psfrag}
\usepackage{pifont}
\usepackage{bbding}
\usepackage[pagebackref]{hyperref}

\definecolor{myred}{rgb}{1,0,0}

\newtheorem{theorem}{Theorem}

{\theoremstyle{remark} \newtheorem{remark}{Remark}}

\newcommand\Eqn{Equation }

\newcommand\Fig{Figure }

\newcommand{\vect}[1]{\textrm{\boldmath${#1}$}} 
\newcommand{\Vect}[1]{{\bf #1}} 

\graphicspath{{images/}}

\newcommand\showfigures[1]{}

\makeatletter
\newenvironment{subfigures}
 {\begin{minipage}{\columnwidth}\def\@captype{figure}\centering}
 {\end{minipage}}
\makeatother

\title{Projective Integration Methods in the Runge-Kutta Framework and the Extension to Adaptivity in Time}
\author{
Julian Koellermeier\footnote{Corresponding author, email address {\tt j.koellermeier@rug.nl}} \footnote{Bernoulli Institute, University of Groningen},
Giovanni Samaey\footnote{Department of Computer Science, KU Leuven}
}

\begin{document}
\maketitle


\begin{abstract}
    Projective Integration methods are explicit time integration schemes for stiff ODEs with large spectral gaps. In this paper, we show that all existing Projective Integration methods can be written as Runge-Kutta methods with an extended Butcher tableau including many stages. We prove consistency and order conditions of the Projective Integration methods using the Runge-Kutta framework. Spatially adaptive Projective Integration methods are included via partitioned Runge-Kutta methods. New time adaptive Projective Integration schemes are derived via embedded Runge-Kutta methods and step size variation while their accuracy, stability, convergence, and error estimators are investigated numerically.
\end{abstract}

\medskip\noindent

\section{Introduction}
Many problems in science and engineering include the solution of stiff ODEs \cite{Hairer2013} of the following form
\begin{equation}
\label{e:ODE}
    \frac{d u}{d t} = f(u,t,\epsilon),
\end{equation}
where $t \in [t_0,T]$ is the time variable, $u \in \mathbb{R}^n$ is the unknown solution, and the right hand side function $f \in \mathbb{R}^n$ is sufficiently smooth and scales with $\mathcal{O}(f) = \frac{1}{\epsilon}$, where the stiffness is governed by a parameter $0 < \epsilon \ll 1$. Small values of $\epsilon$ thus lead to a stiff problem \eqref{e:ODE} with fast modes requiring small time steps of explicit time integration schemes. While the solution of \eqref{e:ODE} is a function of time, note that the right-hand side function might originate from the discretization of spatial derivatives, in which case the left-hand side would be a partial derivative. For the rest of this paper, we assume autonomous ODEs and write $f(u,t,\epsilon)=f(u)$ for conciseness.

Particularly difficult test cases include a scale separation with fast and slow modes present in the solution $u$, occurring in many physical applications \cite{Amrita2022,Hairer2013,Koellermeier2021}. This is the case if the eigenvalues of the right-hand side derivative $\frac{\partial f}{\partial u}$ are given by two clusters around $0$ and $-\frac{1}{\epsilon}$. Even though the slow modes (represented by the modes clustered around $0$) might be most interesting, e.g., in fluid dynamics \cite{Amrita2022,Koellermeier2021,Lafitte2016}, the existing fast modes (represented by the modes clustered around $-\frac{1}{\epsilon}$) need to be resolved in a stable way. Explicit time stepping methods thus result in a severe time step constraint for stable integration of the fast modes. Implicit time stepping methods require a possibly expensive non-linear solver.

Runge-Kutta (RK) methods are well-known methods for the solution of ODEs. Standard explicit RK methods are not suitable to solve stiff ODEs, due to the time step constraint. Modified RK methods like Implicit-Explicit (IMEX) RK methods, are used if the right hand side function of the ODE can be split into a stiff and a non-stiff term, which can be solved separately in an efficient way \cite{Pareschi2005}. A similar splitting is necessary to use Multirate Runge-Kutta (MRK) methods \cite{Gunther2016}. Exponential Runge-Kutta (ERK) methods can be applied for stiff ODEs but have difficulties with scale separation \cite{Hochbruck2020}.

Projective Integration (PI) is an explicit numerical method for the time integration of stiff ODEs that include a scale separation \cite{Gear2003,Lafitte2010,Melis2017}. Its main idea is to first damp the fast modes using a few small time steps of an inner integrator and then extrapolate using a large time step for the long time behavior. For stability, the small (inner) time step size $\delta t$ is of the order of the stiffness parameter $\epsilon \ll 1$ and the large (outer) time step size $\Delta t$ is typically chosen according to a standard CFL condition that corresponds to an underlying limiting equation for $\epsilon \rightarrow 0$, which might not be explicitly given.

PI schemes have been used successfully for many applications, including multiscale ODEs \cite{Maclean2015}, power systems \cite{Wang2018}, stochastic DNA models \cite{Fahrenkopf2013}, shallow water flows \cite{Amrita2022}, and kinetic equations \cite{Lafitte2010,Melis2017,Koellermeier2021}. The method has been extended to higher-order using RK methods for the outer integrator \cite{Lee2007a,Lafitte2017}, which resulted in Projective Runge-Kutta (PRK) methods. Using more projective levels results in the Telescopic Projective Integration (TPI) method for multiple eigenvalue clusters \cite{Melis2016,Melis2019}. A spatially adaptive version was first derived in \cite{KoellermeierAPI}.

In this paper we show that all of the aforementioned PI methods can be written as RK methods and derive the corresponding Butcher tableaus. The inner steps and the subsequent extrapolation step translate to many stages of the corresponding RK method and this generalizes to the higher-order, telescopic, and spatially adaptive versions of the PI method. Writing PI methods in the framework of standard RK schemes allows to make use of the vast number of mathematical tools developed for RK methods throughout the past decades. This includes a simple assessment of consistency and order conditions up to higher order as well as accuracy and numerical stability. Using the RK setup we derive time adaptive PI schemes with the help of embedded RK methods. We compare the new embedded PI methods with other time adaptive methods that are obtained via step size variation, which is a variant of Richardson extrapolation, and on-the-fly error estimation \cite{Gear2006}. Interestingly, we show that the original on-the-fly error estimation from \cite{Gear2006} leads to an instable scheme and we derive a new versions of the schemes using the RK framework for which we prove stability for stiff systems with spectral gaps.
Lastly, we apply the rewritten and newly derived methods do a two-scale model problem and investigate convergence and error estimators to confirm the analytical results obtained before.

While this paper focuses mostly on the method derivation and analysis, it is the necessary preparation for further extensions, more detailed stability analysis, and additional numerical applications, for example using the tools for stability analysis of RK schemes from \cite{Ketcheson2012} or test cases from fluid dynamics.

The rest of this paper is organized as followed: In Section \ref{sec:PI_as_RK}, we briefly recall the general form of a RK method before deriving the corresponding RK formulas for standard PI methods: the first order Projective Forward Euler (PFE) method, the higher order Projective Runge-Kutta (PRK) method, and the telescopic Projective Forward Euler (TPFE) method. Spatially Adaptive Projective Integration (SAPI) methods are written in the framework of partitioned RK methods in Section \ref{sec:SAPI}. Time adaptive Projective Integration (TAPI) methods are derived based on embedded RK methods and step size variation using error estimators in Section \ref{sec:TAPI}. We use the techniques developed in the first section to rewrite the on-the-fly error estimation schemes as a RK scheme and develop a new version of the scheme in Section \ref{on-the-fly}. We perform accuracy and stability analysis in Section \ref{sec:analysis} and numerically assess the convergence properties as well as the error estimators in Section \ref{sec:conv_test}. The paper ends with a short conclusion and future work.

\section{Projective Integration as Runge-Kutta}
\label{sec:PI_as_RK}
In this section, we write recall standard Runge-Kutta (RK) schemes and Projective Integration (PI) methods and write various PI methods as Runge-Kutta schemes including the corresponding Butcher tableau: the first order Projective Forward Euler (PFE) method, the higher order Projective Runge-Kutta (PRK) method, and the telescopic Projective Forward Euler (TPFE) method.

\subsection{Runge-Kutta schemes}
\label{sec:RK}
A standard Runge-Kutta scheme with $S\in \mathbb{N}$ stages computes the intermediate function evaluations $k_j$, for $j=1,\ldots,S$ at every stage according to
\begin{equation}\label{e:RK_stage}
    k_j = f\left(t_n + \Delta t c_j, u + \Delta t \sum_{l=1}^{S} a_{j,l} k_l\right)
\end{equation}
and combines those to a new time iterate
\begin{equation}\label{e:RK_update}
    u^{n+1} = u^{n} + \Delta t \sum_{j=1}^{S} b_j k_j.
\end{equation}
In this setting, the $c_j$ are the intermediate nodes, the $b_j$ are the weights, and the $a_{j,l}$ are the coefficients for the computation of the intermediate function evaluations $k_j$.
Runge-Kutta schemes are typically represented in a Butcher tableau
\begin{equation}
\begin{array}{ c | c}
 \vect{c} & \vect{A} \\
 \hline
  & \vect{b}
  \label{tab:Butcher_RK}
\end{array}
\end{equation}
with nodes vector $\vect{c} \in \mathbb{R}^S$, weight vector $\vect{b} \in \mathbb{R}^S$ and coefficient matrix $\vect{A} \in \mathbb{R}^{S \times S}$. In this paper, we will exclusively consider explicit schemes for which $\vect{A}$ is a lower triangular matrix. Additionally, we note the following conditions for standard properties of explicit Runge-Kutta schemes:
\begin{itemize}
  \item consistency: $\displaystyle\sum_{j=1}^{i-1} a_{i,j} = c_i \quad \textrm{for } i=1,\ldots,S $
  \item first-order accuracy: $\displaystyle\sum_{j=1}^{S} b_{j} = 1 $
  \item second-order accuracy: $\displaystyle\sum_{j=1}^{S} b_{j}c_j = \frac{1}{2} $
\end{itemize}

\subsection{Projective Forward Euler as Runge-Kutta}
\label{sec:PFE_as_RK}
The simplest form of a PI method is the Projective Forward Euler method (PFE) \cite{Gear2003,Melis2017}. It consists of $K+1$ small forward Euler time steps of size $\delta t$ and one subsequent extrapolation step of the remaining time step $\Delta t - (K+1)\delta t$. The method is written using $u^{n,0}=u^n$ and
\begin{eqnarray}
  u^{n,1} &=& u^{n,0} + \delta t \cdot f\left(u^{n,0}\right) \\
  u^{n,2} &=& u^{n,1} + \delta t \cdot f\left(u^{n,1}\right) \\
   &\vdots&  \\
  u^{n,K+1} &=& u^{n,K} + \delta t \cdot f\left(u^{n,K}\right) \\
  u^{n+1} &=& u^{n,K+1} + (\Delta t - (K+1)\delta t) \cdot \frac{u^{n,K+1}-u^{n,K}}{\delta t}.
\end{eqnarray}

To write the PFE in the standard form of a RK scheme, we use the following notation
\begin{eqnarray}
  k_0 &=& f\left(u^{n,0}\right) \\
  k_1 &=& f\left(u^{n,1}\right) = f\left(u^{n,0} + \delta t \cdot k_0 \right) \\
  k_2 &=& f\left(u^{n,2}\right) = f\left(u^{n,0} + \delta t \cdot (k_0 + k_1) \right) \\
   &\vdots&  \\
  k_K &=& f\left(u^{n,K}\right) = f\left(u^{n,0} + \delta t \cdot \sum_{l=0}^{K-1} k_l \right),
\end{eqnarray}
so that the stages $k_l$, for $l=0,\ldots,K$ are the slopes or the right-hand side evaluations at the intermediate points, i.e.
\begin{equation}
  k_l = \frac{u^{n,l+1}-u^{n,l}}{\delta t}.
\end{equation}

We note that the extrapolation step in fact reads
\begin{equation}\label{eq:extrapolation_deltat}
  u^{n+1} = u^{n,K} + (\Delta t - K\delta t) \cdot k_K.
\end{equation}
While this looks significantly simpler than the previous formula for the extrapolation step, still $K+1$ function evaluations for $k_j$ with $j=0,\ldots, K$ are necessary.
In a similar way, we can write the extrapolation step as follows
\begin{eqnarray}
  u^{n+1} &=& u^{n,K+1} + (\Delta t - (K+1)\delta t) \cdot \frac{u^{n,K+1}-u^{n,K}}{\delta t} \\
          &=& u^{n} + \delta t \sum_{l=0}^{K} k_l + (\Delta t - (K+1)\delta t) \cdot k_K \\
          &=& u^{n} + \Delta t \left( \sum_{l=0}^{K-1} \frac{\delta t}{\Delta t} k_l + \left(1 - K \frac{\delta t}{\Delta t}\right) \cdot k_K \right).
\end{eqnarray}

We use the notation $\lambda := \frac{\delta t}{\Delta t}$ in the following. Note that stiff problems $\epsilon \ll 1$ typically require $\delta t \ll \Delta t$, such that $\lambda \ll 1$.

The Butcher tableau of the corresponding RK scheme then reads
\begin{equation}
\begin{array}{ c | c c c c }
 0 & 0 & & & \\
 1 \cdot \lambda & \lambda & 0 & & \\
 \vdots & \vdots & \ddots & \ddots & \\
 K \cdot \lambda & \lambda & \dots & \lambda & 0 \\
 \hline
  & \lambda & \dots & \lambda & 1 - K \lambda
  \label{tab:Butcher_PFE_as_RK}
\end{array}
\end{equation}

For a concise notation of the Butcher tableau in the next sections, we use the following definitions
\begin{equation}
    \vect{\lambda} = \lambda \left(0,1,\ldots,K\right)^T \in \mathbb{R}^{K+1},
\end{equation}
\begin{equation}
    \vect{A_{\lambda}} = {\lambda} \left(\begin{matrix}
                         0 & & & \\
                         1 & 0 & & \\
                         \vdots & \ddots & \ddots & \\
                         1 & \dots & 1 & 0
                       \end{matrix}\right) \in \mathbb{R}^{K+1\times K+1},
\end{equation}
\begin{equation}
    \vect{b_\lambda} = \lambda \left(1,\ldots,1,\frac{1}{\lambda}-K\right)^T \in \mathbb{R}^{K+1},
\end{equation}
such that the Butcher tableau can be written in concise form as
\begin{equation}
\begin{array}{ c | c}
 \vect{\lambda} & \vect{A_\lambda} \\
 \hline
  & \vect{b_\lambda}
  \label{tab:Butcher_PFE}
\end{array}
\end{equation}

For consistency and first order accuracy of the PFE schemes, we refer to the later theorem \ref{th:consistency_PRK}, which proves this for PRK schemes as a superset of PFE methods.

\subsection{Projective Runge-Kutta as Runge-Kutta}
\label{sec:PRK}
Originating from the first-order accurate PFE, higher-order schemes can be constructed using a standard $S$ stage RK method as outer integrator with usual parameters $\vect{A} \in \mathbb{R}^{S \times S}$, $\vect{c} \in \mathbb{R}^{S}$, and $\vect{b} \in \mathbb{R}^{S}$. The result is a Projective Runge-Kutta scheme (PRK) \cite{Lee2007a,Lafitte2016}, consisting of $S$ outer stages, which each include $K+1$ small inner time steps of size $\delta t$ and one subsequent extrapolation step of the remaining time step $c_s \Delta t - (K+1)\delta t$.

The first stage's slopes are computed as
\begin{equation}
    s=1: \left\{
                \begin{array}{lcl}
                  u^{n,k+1} &=& u^{n,k} + \delta t f\left(u^{n,k}\right), \quad 0 \leq k \leq K, \\
                  k_1 &=& \frac{u^{n,K+1}-u^{n,K}}{\delta t}
                \end{array}
              \right.
\end{equation}
and the other stages are subsequently given by
\begin{equation}
    2 \leq s \leq S: \left\{
                \begin{array}{lcl}
                  f_s^{n+c_s,0} &=& u^{n,K+1} + \left( c_s \Delta t - (K+1)\delta t \right) \displaystyle\sum_{l=1}^{s-1} \frac{a_{s,l}}{c_s} k_l \\
                  f_s^{n+c_s,k+1} &=& f_s^{n+c_s,k} + \delta t f\left(f_s^{n+c_s,k}\right), \quad 0 \leq k \leq K \\
                  k_s &=& \frac{f_s^{n+c_s,K+1}-f_s^{n+c_s,K}}{\delta t}
                \end{array}
              \right.
\end{equation}
The new time step is then extrapolated to
\begin{equation}
    u^{n+1} = u^{n,K+1} + (\Delta t - (K+1)\delta t) \sum_{s=1}^{S} b_s k_s.
\end{equation}

To write the PRK in the standard RK form, we use the following notation
\begin{eqnarray}
  s=1 : k_{1,k+1} &=& f\hspace{-0.10cm}\left(\hspace{-0.15cm}u^{n,0} + \delta t \sum_{l=1}^{k} k_{1,l}\right)\\
  2 \leq s \leq S : k_{s,k+1} &=& f\hspace{-0.10cm}\left(\hspace{-0.15cm}u^{n,0} + \delta t \sum_{l=1}^{K+1} k_{1,l} + \left( c_s \Delta t - (K+1)\delta t \right) \sum_{i=1}^{s-1} \frac{a_{s,i}}{c_s} k_{i,k+1}\right.\\
  &\hspace{-0.3cm}\hspace{-0.3cm}& \hspace{1cm} \left. + \delta t \sum_{l=1}^{k} k_{s,l} \right),
\end{eqnarray}
for inner iterations $0 \leq k \leq K$. The new time step is given by
\begin{equation}
    u^{n+1} = u^{n,0} + \Delta t \left( \frac{\delta t}{\Delta t} \sum_{l=1}^{K+1} k_{1,l} + \left(1 - (K+1)\frac{\delta t}{\Delta t}\right) \sum_{s=1}^{S} b_s k_{s,K+1} \right).
\end{equation}

Adapting the notation of the PFE method \ref{sec:PFE_as_RK}, the PRK method can then be written as RK method in a Butcher tableau in the following way
\begin{equation}
    \begin{array}{ c | c c c c c }
     \vect{\lambda} & \vect{A_\lambda} & & & & \\
     \vect{c_2 + \lambda} & \vect{\Lambda^a_{2,1}} & \vect{A_\lambda} & & & \\
     \vect{c_3 + \lambda} & \vect{\Lambda^a_{3,1}} & \vect{0^a_{3,2}} & \vect{A_\lambda} & & \\
     \vdots & \vdots & \vdots & \ddots & \ddots & \\
     \vect{c_S + \lambda} & \vect{\Lambda^a_{S,1}} & \vect{0^a_{S,2}} & \dots & \vect{0^a_{S,S}} & \vect{A_\lambda} \\
     \hline
      & \vect{b_1} & \vect{b_2} & \dots & \vect{b_{S-1}} & \vect{b_S}
    \end{array}
    =:
    \begin{array}{ c | c }
     \overline{\vect{c}} & \overline{\vect{A_\lambda}} \\
     \hline
      & \overline{\vect{b}}
    \end{array}
    \label{tab:Butcher_PRK_as RK}
\end{equation}
with
\begin{equation}
    \vect{c_s} = c_s \left(1,\ldots,1\right)^T \in \mathbb{R}^{K+1}
\end{equation}
\begin{equation}
    \vect{\Lambda^a_{s,1}} = \lambda \left(\begin{matrix}
                         1 & \dots & 1 & 1 + \widetilde{a}_{s,1}\\
                         1 & \dots & 1 & 1 + \widetilde{a}_{s,1}\\
                         \vdots & \ddots & \vdots & \vdots\\
                         1 & \dots & 1 & 1 + \widetilde{a}_{s,1}\\
                       \end{matrix}\right) \in \mathbb{R}^{K+1\times K+1}
\end{equation}
\begin{equation}
    \vect{0^a_{s,l}} = \lambda \left(\begin{matrix}
                         0 & \dots & 0 & \widetilde{a}_{s,l}\\
                         0 & \dots & 0 & \widetilde{a}_{s,l}\\
                         \vdots & \ddots & \vdots & \vdots\\
                         0 & \dots & 0 & \widetilde{a}_{s,l}\\
                       \end{matrix}\right) \in \mathbb{R}^{K+1\times K+1}
\end{equation}
where the last column contains the entry
\begin{equation}
    \widetilde{a}_{s,l} = \left(\frac{c_s}{\lambda} - (K+1)\right) \frac{a_{s,l}}{c_s},
\end{equation}
and the weights are given by
\begin{equation}
    \vect{b_1} = \lambda \left(1,\ldots,1, 1 + \left(\frac{1}{\lambda} - (K+1) \right) b_1 \right)^T \in \mathbb{R}^{K+1},
\end{equation}
\begin{equation}
    \vect{b_s} = \lambda \left(0,\ldots,0, \left(\frac{1}{\lambda} - (K+1) \right) b_s \right)^T \in \mathbb{R}^{K+1}, \quad 2 \leq s \leq S.
\end{equation}

Consistency and order conditions of a PRK method can be easily derived using the representation as a RK method. It is one example for the advantage of writing the existing PI schemes as RK schemes, since the rewriting allows to use the consistency conditions of RK schemes for PI schemes. This is shown in the next theorem.
\begin{theorem}(Consistency and order conditions)
    A PRK method is consistent and at least first order accurate.
    The method recovers second order accuracy of the outer RK method in the case of vanishing $\lambda= \frac{\delta t}{\Delta t} \rightarrow 0$.
    \label{th:consistency_PRK}
\end{theorem}
\begin{proof}
    We proof Theorem \ref{th:consistency_PRK} up to second order by showing that the RK version of the method fulfills the consistency and order conditions.

    Consistency demands $\sum_{j=0}^{S\cdot(K+1)} \overline{a}_{i,j} \overset{!}{=} \overline{c}_i$. We check the condition for each outer stage. For the first outer stages $s=1$, i.e. $i=0,\ldots,K+1$, we obtain for each inner stage $k=0,\ldots,K+1$
    \begin{eqnarray*}
    \sum_{j=0}^{S\cdot(K+1)} \overline{a}_{k,j} &=& \overline{c}_k, \\
      k \cdot \lambda &\overset{!}{=}& \lambda + \ldots + \lambda = k \cdot \lambda \quad\checkmark
    \end{eqnarray*}
    For the later outer stages $s=2,\ldots,S$, we obtain for each inner stage $k=0,\ldots,K+1$
    \begin{eqnarray*}
    \sum_{j=0}^{S\cdot(K+1)} \overline{a}_{(s-1)(K+1)+k,j} &=& \overline{c}_{(s-1)(K+1)+k}, \\ 
      c_s + k \lambda       &=& \underbrace{\lambda + \ldots + \lambda}_{K+1 \text{ times}} + \lambda \sum_{i=1}^{S-1} \widetilde{a}_{s,l} + k\lambda \\
                            &=& (K+1)\lambda + \lambda \sum_{i=1}^{S-1} \left(\frac{c_s}{\lambda} - (K+1)\right) \frac{a_{s,l}}{c_s} + k\lambda \\
                            &=& (K+1)\lambda + \left(c_s - (K+1)\lambda\right) \underbrace{\sum_{i=1}^{S-1} \frac{a_{s,l}}{c_s}}_{=1} + k\lambda \\
                            &=& (K+1)\lambda + \left(c_s - (K+1)\lambda\right) \cdot 1 + k\lambda = c_s + k\lambda, \quad\checkmark
    \end{eqnarray*}
    where we have used that the outer RK scheme is consistent and fulfills $\sum_{i=1}^{S-1} \frac{a_{s,l}}{c_s}=1$.

    First order accuracy demands $\sum_{j=0}^{S\cdot(K+1)} \overline{b}_{j} \overset{!}{=} 1$ and this can be checked by
    \begin{eqnarray*}
    \sum_{j=0}^{S\cdot(K+1)} \overline{b}_{j} &=& \lambda \left( \underbrace{1 + \ldots + 1}_{K+1 \text{ times}}  + \left(\frac{1}{\lambda} - (K+1)\right) b_1 \right) + \sum_{i=2}^{S}\lambda \left( 0 + \ldots + 0  + \left(\frac{1}{\lambda} - (K+1)\right) b_1 \right) \\
    &=& (K+1)\lambda + \sum_{i=1}^{S} (1-(K+1)\lambda)b_i = (K+1)\lambda + (1-(K+1)\lambda) \underbrace{\sum_{i=1}^{S} b_i}_{=1} \\
    &=& (K+1)\lambda + (1-(K+1)\lambda) = 1, \quad\checkmark
    \end{eqnarray*}
    where we used that the outer RK scheme is first order accurate and fulfills $\sum_{i=1}^{S} b_i=1$.

    Second order accuracy demands $\sum_{j=1}^{S\cdot(K+1)} \overline{b}_{j} \overline{c}_{j} \overset{!}{=} \frac{1}{2}$ and we compute
    \begin{eqnarray*}
    \sum_{j=1}^{S\cdot(K+1)} \overline{b}_{j} \overline{c}_{j} &=& \sum_{j=1}^{S} \vect{b_j} (\vect{\lambda}+\vect{c_j}) \\
    &=& \sum_{i=0}^{K} i \lambda^2 + (1-(K+1)\lambda)b_1 K\lambda + \sum_{j=2}^{S} (1-(K+1)\lambda)b_j (c_j + K\lambda) \\
    &=& \lambda^2 \frac{K^2+K}{2} +  \underbrace{\sum_{i=1}^{S} b_i}_{=1} (1-(K+1)\lambda) K\lambda + \underbrace{\sum_{i=1}^{S} b_i c_i}_{=1/2} (1-(K+1)\lambda) \\
    &=& \lambda^2 \frac{K^2+K}{2} + (1-(K+1)\lambda) K\lambda + \frac{1}{2} (1-(K+1)\lambda) = \frac{1}{2} + \mathcal{O}(\lambda).
    \end{eqnarray*}
    Where we have used that the outer RK scheme is second order accurate and therefore fulfills $\sum_{i=1}^{S} b_i=1$ as well as  $\sum_{i=1}^{S} b_i c_i= \frac{1}{2}$.

    The scheme is formally only second order if $\lambda= \frac{\delta t}{\Delta t} \rightarrow 0$. This is due to the first order inner integrator. However, since $\delta t \ll \Delta t$ by potentially several orders of magnitude, the error of the inner integrator can often be neglected. In the limit $\lambda= \frac{\delta t}{\Delta t} \rightarrow 0$ we obtain full second order.
\end{proof}

We note that due to the first order accuracy of the inner integrator, no higher order can formally be obtained in the PRK method. However, the method can yield higher order in the limit of vanishing $\lambda= \frac{\delta t}{\Delta t} \rightarrow 0$ if higher order inner integrators are used.

From the Butcher tableau representation, it can be shown that the original outer RK scheme is recovered (with some redundant zero rows) in the limiting case of vanishing $\lambda= \frac{\delta t}{\Delta t} \rightarrow 0$.

\begin{theorem}(PRK limit)
    In the limit of vanishing $\lambda= \frac{\delta t}{\Delta t} \rightarrow 0$ and constant inner steps $K$, the PRK scheme reduces to its outer RK scheme.
    \label{th:PRK_limit}
\end{theorem}
\begin{proof}
    We show that in the limit of $\lambda= \frac{\delta t}{\Delta t} \rightarrow 0$, the PRK Butcher tableau collapses to the underlying outer RK Butcher tableau    \begin{equation}
    \begin{array}{ c | c }
        \overline{\vect{c}} & \overline{\vect{A_\lambda}} \\
        \hline
        & \overline{\vect{b}}
    \end{array}
        \Rightarrow
    \begin{array}{ c | c }
        \vect{c} & \vect{A} \\
        \hline
        & \vect{b}
    \end{array}
    \end{equation}

    First we show that all submatrices in $\overline{\vect{A_\lambda}}$ have rank one, because of
    \begin{equation}
        \lambda \widetilde{a}_{s,l} = \left(c_s - (K+1)\lambda\right) \frac{a_{s,l}}{c_s} \rightarrow c_s \frac{a_{s,l}}{c_s} = a_{s,l}.
    \end{equation}
    Furthermore, we obtain
    \begin{equation}
        \vect{b_1} = \left(\lambda,\ldots,\lambda, \lambda + \left(1 - (K+1)\lambda \right) b_1 \right)^T \rightarrow \left(0,\ldots,0, b_1\right)^T,
    \end{equation}
    \begin{equation}
        \vect{b_s} =  \left(0,\ldots,0, \left(1 - (K+1)\lambda \right) b_s \right)^T \rightarrow \left(0,\ldots,0, b_s\right)^T, \quad 2 \leq s \leq S
    \end{equation}
    and
    \begin{equation}
        \vect{c_s + \lambda} \rightarrow \left(c_s,\ldots,c_s \right)^T \in \mathbb{R}^{K+1}, \quad 1 \leq s \leq S
    \end{equation}
    Therefore, the inner stages can be eliminated and the Butcher tableau of the PRK scheme collapses to
    \begin{equation*}
    \begin{array}{ c | c}
     \vect{c} & \vect{A} \\
     \hline
      & \vect{b}
    \end{array}
    \end{equation*}
\end{proof}

\begin{remark}(\textit{Multirate Runge-Kutta methods})
    If $\frac{d u}{d t} = f(u,t,\epsilon)$ allows for a splitting into a fast/stiff part and a slow/non-stiff part, e.g. $\frac{d u}{d t} = f_{\textrm{fast}}(u,t,\epsilon) + f_{\textrm{slow}}(u,t) $, tailored RK methods can be applied that can be written as a multirate generalized additive Runge-Kutta (mGARK) method with an extended Butcher tableau \cite{Gunther2016}. For an autonomous system, the Butcher tableau reads
    \begin{equation}
    \begin{array}{ c | c}
     \vect{A}^{f,f} & \vect{A}^{f,s} \\
     \hline
     \vect{A}^{s,f} & \vect{A}^{s,s} \\
     \hline
     \vect{b}^{f} & \vect{b}^{s}
      \label{tab:Butcher_mGARK}
    \end{array}
    \end{equation}
    The work of this section can readily be extended to mGARK methods in which the stiff integrator uses a PI scheme in its RK form. We do not pursue this further as the main advantage of PI schemes is to be applicable without explicit splitting into a stiff and a non-stiff part. However, it might be an interesting option to investigate systems with multiple relaxation rates for which some can be resolved by a splitting approach and others are dealt with non-intrusively using PI.
\end{remark}

In the appendix \ref{app:PRK4}, we give two explicit examples of fourth order PRK schemes that can easily be derived using the notation in this section, for two different numbers of inner time steps. Both methods will be analyzed with respect to stability and numerical convergence in sections \ref{sec:stab_analysis} and \ref{sec:conv_test}.

\subsection{Telescopic Projective Forward Euler as Runge-Kutta}
\label{sec:TPFE}
If the eigenvalue spectrum of the right-hand side function in \eqref{e:ODE} features at least three distinct clusters, Telescopic Projective Integration (TPI) can be used for an acceleration of the method. TPI is a consistent extension of PI in the sense that it uses a nested PI approach \cite{Melis2016}. It consists of several inner integration layers with different time step sizes that need to be chosen to achieve stability.
For simplicity we only focus on the case of two inner integrators and on the Forward Euler scheme (FE) as inner and outer integrator here. The extension towards a Telescopic Projective Runge-Kutta (TPRK) method is then straightforward. In the TPRK scheme a RK scheme is used as the outermost integrator and PFE is used on all the inner levels, where not accuracy but stability is the limiting factor.

Based on \cite{Melis2016} we write the TPFE scheme as follows
\begin{eqnarray}
  k_{1,k+1} &=& f\left(u^{n} + \sum_{l=1}^{k} h_0 k_{1,l}\right) \\
  k_{2,k+1} &=& f\left(u^{n} + \sum_{l=1}^{K_0+1} h_0 k_{1,l} + M_0 h_0 k_{1,K_0+1} + \sum_{l=1}^{k} h_0 k_{2,l} \right) \\
   &\vdots&  \\
  k_{j,k+1} &=& f\left(u^{n} + \sum_{l=1}^{K_0+1} h_0 k_{1,l} + M_0 h_0 k_{1,K_0+1} + \ldots + \sum_{l=1}^{k} h_0 k_{j,l} \right)
\end{eqnarray}
\begin{eqnarray}
  u^{n+1} &=& u^{n} + \sum_{l=1}^{K_0+1} h_0 k_{1,l} + M_0 h_0 k_{1,k_0+1} + \sum_{l=1}^{K_0+1} h_0 k_{2,l} + M_0 h_0 k_{2,k_0+1}\\
   & & + \sum_{l=1}^{K_0+1} h_0 k_{k_1+1,l} + M_0 h_0 k_{k_1+1,k_0+1}  \\
   & & + \sum_{l=1}^{K_0+1} M_1 h_0 k_{k_1+1,l} + M_1 \left( 1+M_0 \right) h_0 k_{k_1+1,k_0+1},
\end{eqnarray}
where the time step size on each layer $i$ is denoted by $h_i$ and the remaining extrapolation (or projective) time step size is $M_i$. Once the time step size of the innermost layer $h_0$ is chosen the others satisfy
\begin{equation}
    h_i = h_0 \prod_{l=0}^{i-1} (M_l + K_l + 1).
    \label{e:TPFE_step_size}
\end{equation}

In the Butcher tableau, the TPFE is written in the following way for the two inner layers $h_0, h_1$ and outer layer $h_2=\Delta t$ using $\lambda_i = \frac{h_i}{\Delta t}$
\begin{equation}
    \begin{array}{ c | c c c c c }
     \vect{\lambda_0} & \vect{A_{\lambda_0}} & & & & \\
     \lambda_1 + \vect{\lambda_0} & \vect{\Lambda_0} & \vect{A_{\lambda_0}} & & & \\
     2\lambda_1 + \vect{\lambda_0} & \vect{\Lambda_0} & \vect{\Lambda_0} & \vect{A_{\lambda_0}} & & \\
     \vdots & \vdots & \vdots & \ddots & \ddots & \\
     K_1 \lambda_1 + \vect{\lambda_0} & \vect{\Lambda_0} & \vect{\Lambda_0} & \dots & \vect{\Lambda_0} & \vect{A_{\lambda_0}} \\
     \hline
      & \vect{b_1} & \vect{b_2} & \dots & \vect{b_{K_1}} & \vect{b_{K_1+1}}
      \label{tab:Butcher_TPFE_as_RK}
    \end{array}
\end{equation}
with the new notation
\begin{equation}
    \vect{\lambda_0} = \lambda_0 \left(0,1,\ldots,K_0\right)^T \in \mathbb{R}^{K_0+1}
\end{equation}
\begin{equation}
    \vect{A_{\lambda_0}} = {\lambda_0} \left(\begin{matrix}
                         0 & & & \\
                         1 & 0 & & \\
                         \vdots & \ddots & \ddots & \\
                         1 & \dots & 1 & 0
                       \end{matrix}\right) \in \mathbb{R}^{K_0+1\times K_0+1}
\end{equation}
\begin{equation}
    \vect{\Lambda_0} = \lambda_0 \left(\begin{matrix}
                         1 & \dots & 1 & 1 + M_0\\
                         1 & \dots & 1 & 1 + M_0\\
                         \vdots & \ddots & \vdots & \vdots\\
                         1 & \dots & 1 & 1 + M_0\\
                       \end{matrix}\right) \in \mathbb{R}^{K_0+1\times K_0+1}
\end{equation}
for weights
\begin{equation}
    \vect{b_j} = \lambda_0 \left(1,\ldots,1, 1 + M_0 \right)^T \in \mathbb{R}^{K_0+1}, \quad 1\leq j \leq K_1,
\end{equation}
\begin{equation}
    \vect{b_{K_1+1}} = \lambda_0  (1+M_1) \left(1,\ldots,1, 1 + M_0 \right)^T \in \mathbb{R}^{K_0+1}.
\end{equation}

\begin{theorem}(Consistency and order conditions)
    The TPFE method is consistent and first order accurate.
    \label{th:consistency_TPFE}
\end{theorem}
\begin{proof}
    Consistency demands $\displaystyle\sum_{j=0}^{(K_0+1)\cdot(K_1+1)} \overline{a}_{i,j} \overset{!}{=} \overline{c}_i$.

    For inner stages $0 \leq k_0 \leq K_0$, and outer stages $0 \leq k_1 \leq K_1$, we thus compute
    \begin{eqnarray*}
        k_1 \lambda_1 + k_0 \lambda_0 &=& \lambda_0 \left(\underbrace{1 + \ldots + 1}_{k_0 \text{ times}}\right) + k_1 \lambda_0 \left(\underbrace{1 + \ldots + 1}_{K_0 +1 \text{ times}} + M_0 \right)\\
        &=& \lambda_0 k_0 + k_1 \underbrace{\lambda_0 \left(K_0+1+M_0\right)}_{\lambda_1} = k_0 \lambda_0 + k_1 \lambda_1, \quad\checkmark
    \end{eqnarray*}
    where we have used that the outer time step ratio is given by $\lambda_1 = \lambda_0 \left(K_0+1+M_0\right)$.

    First order accuracy demands $\displaystyle\sum_{j=0}^{(K_0+1)\cdot(K_1+1)} \overline{b}_{j} \overset{!}{=} 1$ and we analogously compute
    \begin{eqnarray*}
        \sum_{j=0}^{(K_0+1)\cdot(K_1+1)} \overline{b}_{j} &=& \lambda_0 \sum_{j=0}^{K_1} \left(\underbrace{1 + \ldots + 1}_{K_0+1 \text{ times}}+ M_0\right) + \lambda_0 (1+M_1) \left(\underbrace{1 + \ldots + 1}_{K_0+1 \text{ times}}+ M_0\right)\\
        &=& \lambda_0 K_1 \left(K_0 + 1 + M_0\right) + \lambda_0 (1+M_1) \left(K_0 + 1 + M_0\right)\\
        &=& \lambda_0 \left(K_0 + 1 + M_0\right) \left(K_1 + 1 + M_1\right) = 1, \quad\checkmark
    \end{eqnarray*}
    where we have used that $\lambda_0 \left(K_0 + 1 + M_0\right) \left(K_1 + 1 + M_1\right) = 1$ due to Equation \eqref{e:TPFE_step_size}.
\end{proof}

Note that the TPFE scheme cannot be expected to be more than first order, because the outer scheme is a forward euler scheme. The extension to Telescopic Projective Runge-Kutta (TPRK) methods is straightforward but left out for conciseness.

\section{Space Adaptive Projective Integration as partitioned Runge-Kutta}
\label{sec:SAPI}
In this section, we consider the extension to spatial adaptivity, outlined in \cite{KoellermeierAPI}. To that extent, we consider the case that \eqref{e:ODE} is the result of a spatial discretization on a grid in physical space. The vector of unknowns is then a composition of variables at distinct grid points in physical space. Spatial adaptivity is beneficial if several scales are involved. This means that some variables, e.g., on a subdomain of the physical grid, relax faster than others. We assume that the model is autonomous and consists of a stiff part for $u_L$ and a non-stiff part for $u_R$
\begin{eqnarray}
  \frac{d }{d t} u_L &=& f_L(u_L,u_R), \\
  \frac{d }{d t} u_R &=& f_R(u_L,u_R).
\end{eqnarray}
A partitioned Runge-Kutta scheme for this model is given by the stages
\begin{eqnarray}
  u_L^i &=& u_L + \Delta t \sum_{j=1}^{s} a_{i,j} f_L(t_{n+c_j},u_L^j,u_R^j), \\
  u_R^i &=& u_R + \Delta t \sum_{j=1}^{s} \widetilde{a_{i,j}} f_R(t_{n+\widetilde{c_j}},u_L^j,u_R^j).
\end{eqnarray}
The next step is then combined using the stage values to
\begin{eqnarray}
  u_L^{n+1} &=& u_L + \Delta t \sum_{i=1}^{s} b_{i} f_L(t_{n+c_i},u_L^i,u_R^i), \\
  u_R^{n+1} &=& u_R + \Delta t \sum_{i=1}^{s} \widetilde{b_{i}} f_R(t_{n+\widetilde{c_i}},u_L^i,u_R^i).
\end{eqnarray}
A partitioned Runge-Kutta scheme is typically written using two Butcher tableaus, one for the stiff and one for the non-stiff region, respectively.
\begin{equation*}
    \begin{array}{ c | c}
     \vect{c} & \vect{A} \\
     \hline
      & \vect{b}
    \end{array}\quad \quad
      {\def\arraystretch{1.4}
    \begin{array}{ c | c}
    \widetilde{\vect{c}} & \widetilde{\vect{A}} \\
    \hline
    & \widetilde{\vect{b}}
    \end{array}
    }
\end{equation*}

\subsection{Space Adaptive Forward Euler as partitioned Runge-Kutta}
\label{sec:SAFE}
\begin{figure}[htb!]
    \centering
    \includegraphics[width=0.75\textwidth]{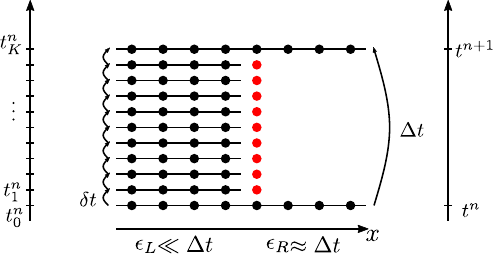}
    \caption{Space Adaptive Forward Euler scheme (SAFE) with small time step $\delta t$ in stiff region (left) and large time step $\Delta t$ in non-stiff region (right). Values of red cells at the boundary of the two domains need to be reconstructed, see \cite{KoellermeierAPI}.}
    \label{fig:AFE_grid}
\end{figure}

As a preparation for the space adaptive PI methods, we first consider the Space Adaptive Forward Euler (SAFE) method from \cite{KoellermeierAPI} as illustrated by Figure \ref{fig:AFE_grid}. The method uses a forward Euler scheme with small step size $\delta t$ for the variables $u_L$ in the stiff domain, and a forward Euler scheme with large time step size $\Delta t$ for the variables $u_R$ in the non-stiff domain. The boundary cells are computed by linear interpolation \cite{KoellermeierAPI}.

This scheme can be written as a partitioned Runge-Kutta method with $M=\frac{\Delta t}{\delta t}$ stages using the two Butcher tableaus
\begin{equation}
    \begin{array}{ c | c c c c c }
     0 & 0 & & & & \\
     1 \cdot \lambda & \lambda & 0 & & & \\
     2 \cdot \lambda & \lambda & \ddots & \ddots & & \\
     \vdots & \vdots &  & \ddots & \ddots &  \\
     M \cdot \lambda & \lambda & \dots & \dots & \lambda & 0 \\
     \hline
      & \lambda & \dots & \dots & \lambda & \lambda
    \end{array}
    \quad\quad
    \begin{array}{ c | c c c c c }
     0 & 0 & & & & \\
     1 \cdot \lambda & 1 \cdot \lambda & 0 & & & \\
     2 \cdot \lambda & 2 \cdot \lambda & \vdots & \ddots & & \\
     \vdots & \vdots & \vdots & & \ddots &  \\
     M \cdot \lambda & M \cdot \lambda & 0 & \dots & \dots & 0 \\
     \hline
      & 1 & 0 & \dots & \dots & 0
    \end{array}
    \label{tab:Butcher_SAFE}
\end{equation}
where the left Butcher tableau in \eqref{tab:Butcher_SAFE} is for the forward Euler scheme with small step size $\delta t$ and the right Butcher tableau in \eqref{tab:Butcher_SAFE} is for the forward Euler scheme with large step size $\Delta t$ including the linear interpolation for the boundary cells.

It is evident from the Butcher tableaus, that the method is first order accurate.

\subsection{Space Adaptive Projective Forward Euler as partitioned Runge-Kutta}
\begin{figure}[htb!]
    \centering
    \includegraphics[width=0.75\textwidth]{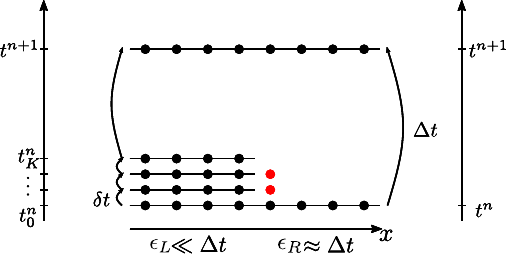}
    \caption{Space Adaptive Projective Forward Euler scheme (SAPFE) with $K+1$ inner small time steps $\delta t$ in stiff region (left) and large time step $\Delta t$ in non-stiff region (right). Values of red cells at the boundary of the two domains need to be reconstructed, see \cite{KoellermeierAPI}.}
    \label{fig:APFE_grid}
\end{figure}

The Space Adaptive Projective Forward Euler (SAPFE) method from \cite{KoellermeierAPI} is illustrated by Figure \ref{fig:APFE_grid}. The method uses a PFE method for the variables $u_L$ in the stiff domain, and Forward Euler method with large time step size $\Delta t$ for the variables $u_R$ in the non-stiff domain. Boundary values for the first inner steps of the PFE method are again obtained by interpolation \cite{KoellermeierAPI}.

Combining the results of Section \ref{sec:PFE_as_RK} and Section \ref{sec:SAFE}, this scheme can be written as a partitioned RK method with $K+1$ stages using the two Butcher tableaus
\begin{equation}
    \begin{array}{ c | c c c c c }
     0 & 0 & & & & \\
     1 \cdot \lambda & \lambda & 0 & & & \\
     2 \cdot \lambda & \lambda & \ddots & \ddots & & \\
     \vdots & \vdots &  & \ddots & \ddots &  \\
     K \cdot \lambda & \lambda & \dots & \dots & \lambda & 0 \\
     \hline
      & \lambda & \dots & \dots & \lambda & 1-K\lambda
        \end{array}
        \quad\quad
        \begin{array}{ c | c c c c c }
     0 & 0 & & & & \\
     1 \cdot \lambda & 1 \cdot \lambda & 0 & & & \\
     2 \cdot \lambda & 2 \cdot \lambda & \vdots & \ddots & & \\
     \vdots & \vdots & \vdots & & \ddots &  \\
     K \cdot \lambda & K \cdot \lambda & 0 & \dots & \dots & 0 \\
     \hline
  & 1 & 0 & \dots & \dots & 0
    \end{array}
    \label{tab:Butcher_SAPFE}
\end{equation}
where the left Butcher tableau in \eqref{tab:Butcher_SAPFE} is for the PFE scheme with inner step size $\delta t$ and outer step size $\Delta t$ and the right Butcher tableau in \eqref{tab:Butcher_SAPFE} is for the forward Euler scheme with large step size $\Delta t$ including the linear interpolation for the boundary cells.

\begin{remark}(\textit{Space Adaptive Projective Runge-Kutta as partitioned RK})
    In the same fashion as the first order PFE method in Section \ref{sec:PFE_as_RK} is extended to higher-order PRK methods in \ref{sec:PRK}, higher-order space adaptive Projective Runge-Kutta schemes can be derived and written in the form of a RK scheme as an extension of the work in this section. However, we omit this here for conciseness and leave specific applications for future work.
\end{remark}

Note also that the consistency and accuracy results can readily be carried over from Section \ref{sec:PRK}.

\subsection{Space Adaptive Projective Projective Forward Euler as partitioned Runge-Kutta}
\begin{figure}[htb!]
    \centering
    \includegraphics[width=0.75\textwidth]{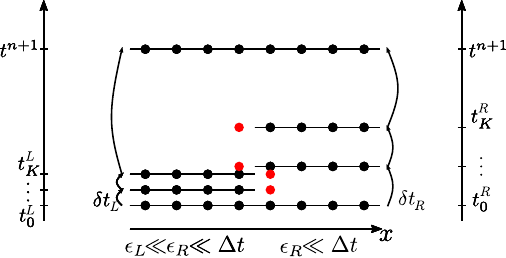}
    \caption{Space Adaptive Projective Projective Forward Euler scheme (SAPPFE) with $K+1=3$ inner small time steps $\delta t_L$ in stiff region (left) \emph{and} $K+1=3$ inner small time steps $\delta t_R > \delta t_L$ in semi stiff region (right). Values of red cells at both sides of the boundary of the two domains need to be reconstructed, see \cite{KoellermeierAPI}. }
    \label{fig:APPFE3_grid}
\end{figure}

The Space Adaptive Projective Projective Forward Euler (SAPPFE) method from \cite{KoellermeierAPI} is illustrated by Figure \ref{fig:APPFE3_grid}. It assumes multiple stiff regions with different severity of time step constraints. The method uses a PFE method in the stiff part of the domain and another PFE method in the semi-stiff part of the domain \cite{KoellermeierAPI}.

This scheme can be written as a partitioned RK method with $(K_L+1)\cdot(K_R+1)$ stages using two Butcher tableaus. The methods first require an ordering of the respective time steps. We use $\lambda_L = \frac{\delta t_L}{\Delta t}$ and $\lambda_R = \frac{\delta t_R}{\Delta t}$ and consider the case $K_L=K_R=2$ as one example. In this case, we obtain:
$$0 < \delta t_L < 2 \delta t_L < \delta t_R < 2 \delta t_R < \Delta t.$$
This fixes the order of the respective stages. The Butcher tableaus then read:
\begin{equation}
    \begin{array}{ c | c c c c c }
     0 & 0 & & & & \\
     \lambda_L & \lambda_L & 0 & & & \\
     2 \lambda_L & \lambda_L & \lambda_L & \ddots & & \\
     \lambda_R & \lambda_L & \lambda_L & \lambda_R - 2\lambda_L & \ddots &  \\
     2 \lambda_R & \lambda_L & \lambda_L & 0 & \lambda_R - 2\lambda_L & 0 \\
     \hline
      & \lambda_L & \lambda_L & 1-2\lambda_L & 0 & 0
      \end{array}
      \quad\quad
      \begin{array}{ c | c c c c c }
     0 & 0 & & & & \\
     \lambda_L & \lambda_L & 0 & & & \\
     2 \lambda_L & \lambda_L & 0 & \ddots & & \\
     \lambda_R & \lambda_R & 0 & 0 & \ddots &  \\
     2 \lambda_R & \lambda_R & 0 & 0 & \lambda_R & 0 \\
     \hline
      & 0 & 0 & \lambda_R & \lambda_R & 1-2\lambda_R
    \end{array}
    \label{tab:Butcher_SAPPFE}
\end{equation}
where the left Butcher tableau in \eqref{tab:Butcher_SAPPFE} is for the PFE scheme with inner step size $\delta t_L$ and outer step size $\Delta t$ including the linear interpolation for the boundary cells and the right Butcher tableau in \eqref{tab:Butcher_SAPPFE} is for the PFE scheme with inner step size $\delta t_R$ and outer step size $\Delta t$ including the linear interpolation for the boundary cells. Note that the right Butcher tableau represents a reducible RK method, due to a full zero row in the matrix $\vect{A}$ and corresponding zero $\vect{b}^T$ entry. However, we note that the values of the right part are needed for the update of the left part and vice versa, see the red dots in Figure \ref{fig:APPFE3_grid}. The method as a whole can therefore not be reduced to less stages.

Consistency and accuracy of this method again follows domain wise from the results in section \ref{sec:PRK}. Note that the extension to more stages or higher-order schemes is straightforward and omitted here for conciseness.

\begin{remark}(\textit{On the equivalence of additive and component partitioning})
    We note that the case of component partitioning that is considered in this section is also included in the case of additive partitioning via a reformulation of the ODE to be solved, see \cite{Gunther2016}. In that sense it is straightforward to extend the results of this section to additively partitioned systems.
\end{remark}

\section{Time Adaptive Projective Integration}
\label{sec:TAPI}
In this section, we derive time adaptive PI methods using three different ways, which are all facilitated by the representation as a RK method derived in section \ref{sec:PI_as_RK}. The novelty of these time adaptive PI methods is that the error estimation of an adaptive time stepping method is based on internal stages of a Projective Integration scheme, facilitated by the formulation as a RK scheme. On the one hand, this allows for a lot of flexibility, e.g., by comparing inner or outer stages with different time step sizes. On the other hand, it allows to redefine existing methods such as embedded RK methods or Richardson extrapolation, as will be discussed in the subsections below.

\subsection{Time Adaptive Projective Integration from embedded Runge-Kutta}
\label{sec:TAPI_embedded}
Embedded RK schemes use a combination of two schemes, one low order and one higher order \cite{Hairer2000}. This leads to an extended Butcher tableau of the following form
\begin{equation}
{\def\arraystretch{1.4}
\begin{array}{ c | c}
 \vect{c} & \vect{A} \\
 \hline
  & \vect{b} \\
 \hline
  & \widetilde{\vect{b}}
  \label{tab:Butcher_eRK}
\end{array}
}
\end{equation}
where the first row containing $\vect{b}$ denotes the higher order RK scheme and the second row containing $\widetilde{\vect{b}}$ denotes the low order RK scheme. Note that PI schemes will only formally be second or higher order for vanishing $\lambda \rightarrow 0$.

As the extended Butcher tableau is nothing else than a combination of two standard RK tableaus with the same $\vect{c}$ and $\vect{A}$, but two different $\vect{b}$ and $\widetilde{\vect{b}}$, the extension to a PRK scheme written in this form is straightforward and reads
\begin{equation}
{\def\arraystretch{1.4}
\begin{array}{ c | c}
 \overline{\vect{c}} & \overline{\vect{A}} \\
 \hline
  & \overline{\vect{b}} \\
 \hline
  & \overline{\widetilde{\vect{b}}}
  \label{tab:Butcher_ePRK}
\end{array}
}
\end{equation}
where the entries $\overline{\vect{c}}$, $\overline{\vect{A}}$, and $\overline{\vect{b}}$ are defined in Section \ref{sec:PRK} and $\overline{\widetilde{\vect{b}}}$ is the corresponding weight vector of the higher order RK method.

Embedded RK schemes allow for a cheap error estimator, as the low order and the high order method use the same nodes $\vect{c}$ with the same coefficients $\vect{A}$ and the difference between the high order solution $u^{n+1}$ and the low order solution $\widetilde{u}^{n+1}$ is given by
\begin{equation}
    u^{n+1} - \widetilde{u}^{n+1} = \Delta t \sum_{j=1}^{s} (b_j - \widetilde{b_j}) k_j,
\end{equation}
where $k_j$ contains the function evaluation at time $t_{n+ c_j}$.

For an embedded PRK method with $S$ outer stages and $K+1$ inner steps, this leads to
\begin{eqnarray}
    u^{n+1} - \widetilde{u}^{n+1} &=& \Delta t \sum_{i=1}^{S(K+1)} \left(\vect{b}_i - \overline{\widetilde{\vect{b}}}_i\right) k_i  \\
    &=& \Delta t \sum_{i=1}^{S} \lambda \left( \frac{1}{\lambda} - (K+1)\right) (b_i - \widetilde{b}_i )k_{(K+1)i}\\
    &=& \left( \Delta t - (K+1) \delta t\right)\sum_{i=1}^{S} (b_i - \widetilde{b}_i ) k_{(K+1)i},
    \label{e:embedded_PRK_estimator}
\end{eqnarray}
which is a consistent extension of the standard embedded RK case combined with the extrapolation step of PI schemes.

As the simplest example of an embedded RK method with only 2 stages, we consider the combination of the Heun method with the Forward Euler method:
\begin{equation}
{\def\arraystretch{1.2}
\begin{array}{ c | c c }
 0 &  & \\
 1 & 1 & \\
 \hline
  & \frac{1}{2} & \frac{1}{2} \\
 \hline
  & 1 & 0
  \label{tab:Butcher_eRK2}
\end{array}
}
\end{equation}
In the example above, the standard error estimator is
\begin{equation}
    u^{n+1} - \widetilde{u}^{n+1} = \Delta t \left(-\frac{1}{2}\right) (k_1 - k_2).
\end{equation}

This embedded method can make use of projective integration by substituting the respective projective versions of the schemes. In the simple case above this means substituting the Projective Heun method for the Heun method and the Projective Forward Euler method for the Forward Euler method. Following the explanation for general PRK schemes in Section \ref{sec:PRK} this leads to the extended Butcher tableau:
\begin{equation}
{\def\arraystretch{1.4}
\begin{array}{ c | c c }
 \vect{\lambda} & \vect{A_\lambda} & \\
 \vect{1 + \lambda} & \vect{\Lambda^a_{2,1}} & \vect{A_\lambda}\\
 \hline
  & \vect{b_1} & \vect{b_2} \\
 \hline
  & \widetilde{\vect{b_1}} & \widetilde{\vect{b_2}}
  \label{tab:Butcher_ePRK2}
\end{array}
}
\end{equation}
where the respective entries are defined in Section \ref{sec:PRK}.

Using $K=2$, the parameters for the Projective Heun method and the Projective Forward Euler method yield the following embedded Runge-Kutta scheme, which we call the \emph{Embedded Projective Heun Projective Forward Euler} scheme (EPHPFE):
\begin{equation}
\begin{array}{ c | c c c c c c }
 0              & 0 & & & & & \\
 \lambda      & \lambda & 0 & & & & \\
 2 \lambda    & \lambda & \lambda & 0 & & & \\
 1 + 0        & \lambda & \lambda & 1 - 2 \lambda & 0 & & \\
 1 + \lambda  & \lambda & \lambda & 1 - 2 \lambda & \lambda & 0 &\\
 1 + 2 \lambda& \lambda & \lambda & 1 - 2 \lambda & \lambda & \lambda & 0 \\
 \hline
                & \lambda & \lambda & \frac{1}{2} - \frac{1}{2} \lambda & 0 & 0 & \frac{1}{2} - \frac{3}{2} \lambda \\
 \hline
                & \lambda & \lambda & 1-2\lambda & 0 & 0 & 0
  \label{tab:Butcher_ePRKK2}
\end{array}
\end{equation}
This is an embedded scheme, which formally results in second order for vanishing $\lambda \rightarrow 0$.

The error estimator is then simply obtained using
\begin{eqnarray}
  u^{n+1} - \widetilde{u}^{n+1} &=& \Delta t \sum_{j=1}^{s} (b_j - \widetilde{b_j}) k_j \\
   &=& \Delta t \left( -\frac{1}{2} + \frac{3}{2} \lambda \right) k_3 + \left( \frac{1}{2} - \frac{3}{2} \lambda \right) k_6 \\
   &=& \Delta t \left( \frac{3}{2} \lambda - \frac{1}{2} \right) (k_3 - k_6)
\end{eqnarray}
and can be seen as an example of Equation \eqref{e:embedded_PRK_estimator}.

The corrected scheme then has the following Butcher tableau
\begin{equation}
\begin{array}{ c | c c c c c c }
 0              & 0 & & & & & \\
 \lambda      & \lambda & 0 & & & & \\
 2 \lambda    & \lambda & \lambda & 0 & & & \\
 1 + 0        & \lambda & \lambda & 1 - 2 \lambda & 0 & & \\
 1 + \lambda  & \lambda & \lambda & 1 - 2 \lambda & \lambda & 0 &\\
 1 + 2 \lambda& \lambda & \lambda & 1 - 2 \lambda & \lambda & \lambda & 0 \\
 \hline
              & \lambda & \lambda & \frac{1}{2} - \frac{1}{2} \lambda & 0 & 0 & \frac{1}{2} - \frac{3}{2} \lambda
  \label{tab:Butcher_ePRKK2_cor}
\end{array}
\end{equation}

The error estimator can also be used to subsequently adjust the next time step based on some control strategy, but this is out of the scope of this work.

\begin{remark}(on nodes $c_i > 1$)
    In a PRK scheme with higher-order outer integrator nodes with $c_i>1$ can occur, see for example \eqref{tab:Butcher_ePRKK2}. This is unusual in view of standard RK schemes. However, the additional stages with $c_i>1$ originate from the inner integrator iterating after a usual outer iteration to damp the fast modes in the solution. This yields the desired stability properties also discussed in section \ref{sec:stab_analysis}. For more details on PRK schemes, we refer to the literature \cite{Lafitte2016,Lafitte2017}.
\end{remark}

\subsection{Time Adaptive Projective Integration from outer step size variation}
\label{sec:POSV}
Another way of deriving an error estimator for adaptive time stepping is via computing the solution again with a smaller step size $\frac{\Delta t}{k}$, for $k \in \mathbb{N}$, also called Richardson extrapolation \cite{Richardson1911}. Assuming a time stepping scheme of order $p$, the solution $u^{n+1}_{\Delta t/k}$ at time $t + \Delta t$ using $k$ steps of size $\frac{\Delta t}{k}$ is deviating from the exact solution $u^{n+1}_{\ast}$ according to
\begin{equation}
    u^{n+1}_{\Delta t/k} = u^{n+1}_{\ast} + k C(t_n) \cdot \left(\frac{\Delta t}{k}\right)^{p+1} + \mathcal{O}\left(\Delta t^{p+2}\right).
\end{equation}
Using two different solutions for $k=1$ and $k=2$, we can eliminate the error constant $C(t_n)$ and obtain the corrected solution
\begin{equation}
    u^{n+1}_{\ast} = \frac{2^p u^{n+1}_{\Delta t/2} - u^{n+1}_{\Delta t/1}}{2^p - 1}.
\end{equation}
In the same way, the error estimator can be computed as
\begin{equation}
    u^{n+1}_{\Delta t/1} - u^{n+1}_{\ast} = \frac{2^p \left(u^{n+1}_{\Delta t/1} - u^{n+1}_{\Delta t/2}\right)}{2^p - 1}.
\end{equation}

As one known example, we consider the case of a Forward Euler scheme with $\Delta t/2$ to obtain $u^{n+1}_{\Delta t/2}$ and with $\Delta t/1$ to obtain $u^{n+1}_{\Delta t/1}$. Writing both methods as a two-stage method, this leads to the following embedded scheme:
\begin{equation}
{\def\arraystretch{1.2}
\begin{array}{ c | c c }
 0 &  & \\
 \frac{1}{2} & \frac{1}{2} & \\
 \hline
  & \frac{1}{2} & \frac{1}{2} \\
 \hline
  & 1 & 0
  \label{tab:Butcher_TAPI_outer_embedded}
\end{array}
}
\end{equation}
for which the error estimate reads
\begin{equation}
    u^{n+1}_{\Delta t/1} - u^{n+1}_{\ast} = 2 \left(u^{n+1}_{\Delta t/1} - u^{n+1}_{\Delta t/2}\right) = \Delta t (k_1 - k_2).
    \label{e:error_estimate_Heun}
\end{equation}

The corrected solution is precisely the solution obtained by the explicit midpoint scheme:
\begin{equation}
\begin{array}{ c | c c }
 0 &  & \\
 \frac{1}{2} & \frac{1}{2} & \\
 \hline
  & 0 & 1
  \label{tab:Butcher_TAPI_outer_corrected}
\end{array}
\end{equation}

\begin{figure}[htb!]
    \centering
    \includegraphics[width=0.5\textwidth]{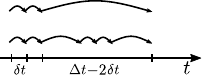}
    \caption{Projective outer step size variation (POSV) for Projective Forward Euler scheme with $K=2$. Top: standard PFE with outer step size $\Delta t$. Bottom: PFE with outer step size $\Delta t/2$.}
    \label{fig:TAPI_outer}
\end{figure}

This scheme can be cast as a projective scheme using the approach for general RK methods from Section \ref{sec:PRK} or for embedded methods from Section \ref{sec:TAPI_embedded}. This includes substituting a Projective Forward Euler method for both Forward Euler methods, once with outer time step $\Delta t$ and once with outer time step $\Delta t/2$. The approach is illustrated in Figure \ref{fig:TAPI_outer} for $K=2$.

The corresponding projective embedded scheme reads
\begin{equation}
\begin{array}{ c | c c c c c c }
 0              & 0 & & & & & \\
 \lambda      & \lambda & 0 & & & & \\
 2 \lambda    & \lambda & \lambda & 0 & & & \\
 \frac{1}{2} + 0        & \lambda & \lambda & \frac{1}{2} - 2\lambda & 0 & & \\
 \frac{1}{2} + \lambda  & \lambda & \lambda & \frac{1}{2} - 2\lambda & \lambda & 0 &\\
 \frac{1}{2} + 2 \lambda& \lambda & \lambda & \frac{1}{2} - 2\lambda & \lambda & \lambda & 0 \\
 \hline
                & \lambda & \lambda & \frac{1}{2} - \frac{1}{2}\lambda & 0 & 0 & \frac{1}{2} - \frac{3}{2}\lambda\\
 \hline
                & \lambda & \lambda & 1 - 2 \lambda & 0 & 0 & 0
  \label{tab:Butcher_TAPI_outer_K2_embedded}
\end{array}
\end{equation}
This is an embedded scheme, which formally results in second order for vanishing $\lambda \rightarrow 0$. The error estimate is
\begin{equation}
    u^{n+1}_{\Delta t/1} - u^{n+1}_{\ast} = \Delta t \left( - \frac{1}{2} + \frac{3}{2}\lambda \right) (k_3 - k_6).
\end{equation}
Note the consistent similarity with the original error estimator in Equation \eqref{e:error_estimate_Heun} and the similarity to the embedded scheme in Section \ref{sec:TAPI_embedded}, which uses the same weights $\overline{\vect{b}}$ and $\overline{\widetilde{\vect{b}}}$.

The Butcher tableau for the new corrected method, which we call \emph{Projective Outer Step Size Variation} (POSV), reads
\begin{equation}
\begin{array}{ c | c c c c c c }
 0              & 0 & & & & & \\
 \lambda      & \lambda & 0 & & & & \\
 2 \lambda    & \lambda & \lambda & 0 & & & \\
 \frac{1}{2} + 0        & \lambda & \lambda & \frac{1}{2} - 2\lambda & 0 & & \\
 \frac{1}{2} + \lambda  & \lambda & \lambda & \frac{1}{2} - 2\lambda & \lambda & 0 &\\
 \frac{1}{2} + 2 \lambda& \lambda & \lambda & \frac{1}{2} - 2\lambda & \lambda & \lambda & 0 \\
 \hline
                & \lambda & \lambda & 0 & 0 & 0 & 1 - 2 \lambda
  \label{tab:Butcher_TAPI_outer_K2_corrected}
\end{array}
\end{equation}
For the stability analysis and numerical convergence of this POSV scheme, see sections \ref{sec:stab_analysis} and \ref{sec:conv_test}.

While a generalization to other PI methods based on different outer RK methods and to arbitrary number of inner steps $K+1$ is possible, we omit this here for conciseness.

\subsection{Time Adaptive Projective Integration from inner step size variation}
\label{sec:PISV}
\begin{figure}[htb!]
    \centering
    \includegraphics[width=0.5\textwidth]{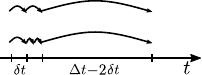}
    \caption{Last projective inner step size variation (PISV) for Projective Forward Euler scheme. Top: standard PFE with inner step size $\delta t$. Bottom: PFE with last inner step using two steps of size $\delta t/2$.}
    \label{fig:TAPI_inner}
\end{figure}

Another possibility to estimate the error of a PI scheme is to change the step size of the inner time step. In principle, all inner time steps could be performed with a smaller inner time step size. However, it is also sufficient to only apply this to the last inner time step. Figure \ref{fig:TAPI_inner} illustrates the idea.

In order to derive an embedded scheme and the corresponding error estimator in the framework of RK methods, we start from the RK version of the single schemes, see Section \ref{sec:PFE_as_RK}. The coarse scheme is a standard PFE method. The fine scheme is an adapted PFE method that uses two smaller inner time steps before the extrapolation step, see Figure \ref{fig:TAPI_inner}. For simplicity, we use $K=1$, but the results can readily be extended to arbitrary $K$. This leads to the following Butcher tableaus for both schemes
\begin{equation}
    \begin{array}{ c | c c }
     0 & & \\
     \lambda & \lambda & \\
     \hline
      & \lambda & 1 - \lambda
    \end{array}\quad \quad
    \begin{array}{ c | c c c }
     0 & & & \\
     \lambda & \lambda & & \\
     1.5 \lambda & \lambda & \frac{1}{2}\lambda & \\
     \hline
      & \lambda & \frac{1}{2}\lambda & 1 - \frac{3}{2} \lambda
    \end{array}
      \label{tab:Butcher_TAPI_inner_two}
\end{equation}
Obviously, the methods do not compute the same stages. This problem can be mitigated by adding a (redundant) additional stage $1.5\lambda$ in the coarse method. This leads to the new embedded scheme:
\begin{equation}
\begin{array}{ c | c c c }
 0 & & & \\
 \lambda & \lambda & & \\
 1.5 \lambda & \lambda & \frac{1}{2}\lambda & \\
 \hline
  & \lambda & \frac{1}{2}\lambda & 1 - \frac{3}{2} \lambda \\
 \hline
  & \lambda & 1- \lambda & 0
  \label{tab:Butcher_TAPI_inner_embedded}
\end{array}
\end{equation}
The error estimator at $t_n + \Delta t$ is obtained using
\begin{equation}
    u^{n+1}_{\delta t/1} - u^{n+1}_{\ast} = \Delta t \left( -1 + \frac{3}{2}\lambda \right) (k_2 - k_3).
\end{equation}

Using Richardson extrapolation, the error at $t_n + 2\delta t$ can be estimated, too. This error estimator can then be used to correct the coarse solution at $t_n + 2\delta t$. In analogy to the Richardson extrapolation in Section \ref{sec:POSV}, we obtain
\begin{equation}
    u^{n+2\delta t}_{\ast} = 2 u^{n+2\delta t}_{\delta t/2} - u^{n+2\delta t} = u^{n+\delta t}_{\delta t/1} + \delta t f\left( u^{n+1.5\delta t} \right).
\end{equation}
Performing the extrapolation step thereafter as follows
\begin{equation}
    u^{n+1} = u^{n+2\delta t}_{\ast} + (\Delta t - 2 \delta t) \frac{u^{n+2\delta t}_{\ast} - u^{n+\delta t}}{\delta t},
\end{equation}
yields the following corrected RK scheme, which we call \emph{Projective Inner Step Size Variation} (PISV),
\begin{equation}
\begin{array}{ c | c c c }
 0 & & & \\
 \lambda & \lambda & & \\
 1.5 \lambda & \lambda & \frac{1}{2}\lambda & \\
 \hline
  & \lambda & 0 & 1 - \lambda
  \label{tab:Butcher_TAPI_inner_one_corrected}
\end{array}
\end{equation}
Note how the solution at time $t_n + 1.5\delta t$ is used as a corrector to improve the solution.
For the stability analysis and numerical convergence of this PISV scheme, we again refer to sections \ref{sec:stab_analysis} and \ref{sec:conv_test} below.

While presenting an example using the first order PFE scheme here, the idea can be extended to higher-order of the outer integrator using the technique explained in Section \ref{sec:PRK}. We leave the application of higher-order schemes for future work.

\begin{remark}(\textit{Space \emph{and} time adaptive PI})
    Using the tools presented in the preceding sections, it is possible to construct space \emph{and} time adaptive PI schemes from embedded, partitioned RK methods. While the time adaptivity can be realized using embedded RK schemes, the spatial adaptivity is dealt with using partitioned RK schemes. Additional necessary ingredients are an efficient error control strategy for the spatial and temporal adaptivity. The combination of spatial and temporal adaptivity is a powerful possibility to construct efficient numerical schemes for problems that exhibit nonlinear or local phenomena, that need to be resolved both in space and in time. We leave this for future work.
\end{remark}

\section{On-the-fly error estimation with Projective Integration}
\label{on-the-fly}
In \cite{Gear2006}, an alternative method for reducing the error of projective integration schemes is proposed. The method uses a so-called on-the-fly error estimation to measure the leading term of the global error after one full projective integration time step, i.e. the accumulated error of the $K+1$ inner steps and the extrapolation step. The method can then be corrected for a more accurate solution. The error of one projective integration step $u_{n+1}$ starting from a correct solution $u_n$ is defined as
\begin{equation}\label{eq:global_error}
    e_{n+1} = u_{n+1} - u\left( t_{n+1} \right).
\end{equation}

The derivation from \cite{Gear2006} has so far not gained a lot of attention in the literature, probably in parts due to the complicated formulation with a number of additional computations. However, writing the underlying PI scheme in its Runge-Kutta version derived in this paper allows to obtain a concise formulation of the scheme that can readily be applied in any Runge-Kutta solver framework,

We here consider a first order PFE scheme to outline the idea. The extension to higher order PRK schemes can be performed using the derivations in \cite{Gear2006}. The aim of the on-the-fly error estimation method is to measure the leading error term at the newly computed solution time $t_{n+1}$ given by
\begin{equation}\label{eq:leading_error}
    e_{n+1} = -\xi \frac{\Delta t^2}{2} u''\left(t_{n+1}\right) + \mathcal{O}\left(\Delta t^3\right),
\end{equation}
where $\xi$ is the unknown leading error coefficient to be estimated and $u''$ is an approximation to the second derivative of the solution.

The on-the-fly error estimation now considers the two ingredients:
\begin{enumerate}
  \item estimate the leading error coefficient $\xi$ for a PFE using $K+1$ inner steps
  \item approximate the second derivative $u''$ at current time with different methods
\end{enumerate}

As for the first ingredient (1), the authors consider the error propagation of the inner steps $j=0,\ldots,K$, based on the error coefficient of a simple forward Euler method, for which the separate error coefficient would be $1$. The accumulated error coefficient after each inner step is then given by
\begin{equation}\label{eq:error_coef_single_step}
    \xi_j = j
\end{equation}

The next step is to derive the error propagation during the extrapolation step. For this, the authors of \cite{Gear2006} use the notation
\begin{equation*}
    u^{n+1} = u^{n,K+1} + M \left( u^{n,(K+1)} - u^{n,K} \right),
\end{equation*}
which can be written in the form of
\begin{equation*}
    u^{n+1} = u^{n,K+1} + (\Delta t - (K+1)\delta t) \cdot \frac{u^{n,K+1}-u^{n,K}}{\delta t}
\end{equation*}
using
\begin{equation}\label{eq:extrapolation_size}
    M=\frac{1}{\lambda} - K - 1
\end{equation}
as extrapolation step size.

According to \cite{Gear2006}, the propagated error coefficient after the extrapolation step at time $s= (M + K + 1) \delta t$ is then given by
\begin{equation}\label{eq:coef_extrapolated}
    \xi^s = (M+1) \xi_{j+1} - M \xi_{j} +  M (M+1),
\end{equation}
which can be evaluated using \eqref{eq:error_coef_single_step} and \eqref{eq:extrapolation_size} to
\begin{eqnarray}
  \xi_s &=& \left( \frac{1}{\lambda} - K \right) \left( K+1 \right) - \left( \frac{1}{\lambda} - K -1 \right) K + \left( \frac{1}{\lambda} - K - 1 \right) \left( \frac{1}{\lambda} - K \right) \\
   &=& \frac{1}{\lambda^2 } - \frac{2K}{\lambda} + K^2 + K.
\end{eqnarray}

Lastly, the error coefficient needs to be corrected for the fact that it was computed in terms of $\delta t$ and not $\Delta t$ by a simple scaling
\begin{equation}\label{eq:error_coef_Scaled}
    \xi = \frac{\xi_s}{s^2} = 1 - 2K \lambda +  (K^2+K) \lambda^2.
\end{equation}
This means that the leading error coefficient can be estimated without any further function evaluations, as soon as the inner and outer time step sizes and the number of inner time steps are known. Note that the leading error coefficient converges to $\xi =1$ for vanishing $\lambda \rightarrow 0$, consistently reproducing the value of the leading error coefficient for a full Forward Euler step.

As for the second ingredient (2), different choices are possible leading to different schemes, which we will describe separately.

\subsection{Outer derivative approximation On-the-fly Projective Forward Euler (OPFE)}
In \cite{Gear2006}, the authors suggest to approximate the derivative in \eqref{eq:leading_error} using the outer stages and \eqref{e:ODE} as
\begin{equation}
    \Delta t^2 u''\left(t_{n+1}\right) \approx \Delta t \left( f \left(u^{n+1},t_{n+1},\epsilon \right) - f \left(u^{n},t_{n},\epsilon \right)\right).
\end{equation}
This approximation can be obtained at the cost of one additional function evaluation of $f$ at time $t_{n+1}$, which will be reflected by one additional Runge-Kutta stage with $c=1$.

The corrected solution $\widetilde{u^{n+1}}$ is then obtained by simple subtraction of the error \eqref{eq:global_error} as
\begin{equation}
    \widetilde{u^{n+1}} = u^{n+1} - e^{n+1}.
\end{equation}

Using a first order PFE with $K+1$ inner stages and one additional stage for the derivative estimation, this leads to the following Butcher tableau for the outer derivative estimation on-the-fly Projective Forward Euler (OPFE) scheme
\begin{equation}
\begin{array}{ c | c c c c c}
 0 & 0 & & & & \\
 1 \cdot \lambda & \lambda & 0 & & &\\
 \vdots & \vdots & \ddots & \ddots & &\\
 K \cdot \lambda & \lambda & \dots & \lambda & 0 &\\
 1 & \lambda & \dots & \lambda & 1 - K \lambda & 0 \\
 \hline
  & \lambda - \frac{\xi}{2}& \dots & \lambda & 1 - K \lambda & \frac{\xi}{2}
  \label{tab:Butcher_PFE_on_the_fly}
\end{array}
\end{equation}

Note that in the limit $\lambda \rightarrow 0$, we obtain $\xi = 1$, the inner stages vanish and the scheme given by \eqref{tab:Butcher_PFE_on_the_fly} will degenerate to the second-order Heun scheme, given by the Butcher tableau
\begin{equation}
\begin{array}{ c | c c }
 0 &  & \\
 1 & 1 & \\
 \hline
  & \frac{1}{2} & \frac{1}{2}
  \label{tab:Butcher_Heun}
\end{array}
\end{equation}

To analyse the order of the scheme, we check the order conditions from section \ref{sec:PI_as_RK} and can easily prove second-order accuracy, independent of the outer time step size, inner time step size, and number of inner time steps. This makes the method a powerful tool to achieve second order at the expense of only one function evaluation.

For illustration, we consider the three examples $K+1=1,2,3$.

For $K=1$, we obtain from \eqref{eq:coef_extrapolated} that $\xi = 1 - 2\lambda + 2\lambda^2$ and
the Butcher tableau of the OPFE1 reads
\begin{equation}
\begin{array}{ c | c c c}
 0 & 0 & & \\
 \lambda & \lambda & 0 &\\
 1 & \lambda & 1 - \lambda & 0 \\
 \hline
  & - \lambda^2 + 2 \lambda - \frac{1}{2} & 1 - \lambda & \lambda^2 - \lambda + \frac{1}{2}
  \label{tab:Butcher_PFE_on_the_fly_K1}
\end{array}
\end{equation}

For $K=2$, we obtain from \eqref{eq:coef_extrapolated} that $\xi = 1 -4\lambda + 6\lambda^2$ and
the Butcher tableau of the OPFE2 reads
\begin{equation}
\begin{array}{ c | c c c c}
 0 & 0 & & & \\
 \lambda & \lambda & 0 & &\\
 2 \lambda & \lambda & \lambda & &\\
 1 & \lambda & \lambda & 1 - 2 \lambda & 0 \\
 \hline
  & - 3 \lambda^2 + 3 \lambda - \frac{1}{2} & \lambda & 1 - 2 \lambda & 3 \lambda^2 - 2 \lambda + \frac{1}{2}
  \label{tab:Butcher_PFE_on_the_fly_K2}
\end{array}
\end{equation}

For $K=3$, we obtain from \eqref{eq:coef_extrapolated} that $\xi = 1 - 6\lambda + 12\lambda^2$ and
the Butcher tableau of the OPFE3 reads
\begin{equation}
\begin{array}{ c | c c c c c}
 0 & 0 & & & &\\
 \lambda & \lambda & 0 & & &\\
 2 \lambda & \lambda & \lambda & & &\\
 3 \lambda & \lambda & \lambda & \lambda & &\\
 1 & \lambda & \lambda & \lambda & 1 - 3 \lambda & 0\\
 \hline
  & - 6 \lambda^2 + 4 \lambda - \frac{1}{2} & \lambda & \lambda & 1 - 3 \lambda & 6 \lambda^2 - 3 \lambda + \frac{1}{2}
  \label{tab:Butcher_PFE_on_the_fly_K3}
\end{array}
\end{equation}

By checking the consistency and order conditions in \ref{sec:RK}, it can be shown again that all scheme are second order accurate.

We emphasize that the OPFE method is estimating the derivative in \eqref{eq:leading_error} using function evaluations at the outer step values $u^{n+1}$ and $u^n$. However, the usage of the outer steps does not take into account the dynamics of the fast scales in the system, which can negatively impact the stability region, as will be investigated during the stability analysis in section \ref{sec:stab_analysis}. It is thus desirable to consider alternatives to the suggested procedure from \cite{Gear2006}.

\subsection{Inner derivative approximation On-the-fly Projective Forward Euler (IPFE)}
We will now present a new alternative to the estimation of the derivative term in \eqref{eq:leading_error}, different from using the outer step values $u^{n+1}$ and $u^n$ as suggested by \cite{Gear2006}. The motivation is to perform the derivative estimation not over the longer time scale $\Delta t$, but over the shorter time scale $\delta t$, to make sure that the corrected scheme does not suffer from instability for fast modes.

One example for this is to estimate the derivative using an additional micro step of size $\delta t$, according to
\begin{equation}
    \Delta t^2 u''\left(t_{n+1}\right) \approx \Delta t \left( f \left(u^{n+1}+\delta t \cdot f\left(u^{n+1}, t_{n+1}, \epsilon \right),t_{n+1}+\delta t,\epsilon \right) - f \left(u^{n+1},t_{n+1},\epsilon \right)\right).
\end{equation}
This approximation can be obtained at the cost of two additional function evaluations of $f$ at time $t_{n+1}$ and $t_{n+1}+\delta t$, which will be reflected by two additional Runge-Kutta stages with $c=1$ and $c=1+\lambda$.

Using the estimated derivative and the error, the corrected solution $\widetilde{u^{n+1}}$ is again obtained by simple subtraction of the error \eqref{eq:global_error} as
\begin{equation}
    \widetilde{u^{n+1}} = u^{n+1} - e^{n+1}.
\end{equation}

For a standard PFE with $K+1$ inner stages and the two additional stages for the derivative estimation, this leads to the following Butcher tableau for the inner derivative estimation on-the-fly Projective Forward Euler (IPFE) scheme
\begin{equation}
\begin{array}{ c | c c c c c c}
 0 & 0 & & & & &\\
 1 \cdot \lambda & \lambda & 0 & & & &\\
 \vdots & \vdots & \ddots & \ddots & & &\\
 K \cdot \lambda & \lambda & \dots & \lambda & 0 & &\\
 1 & \lambda & \dots & \lambda & 1 - K \lambda & 0 &\\
 1+\lambda & \lambda & \dots & \lambda & 1 - K \lambda & \lambda & 0 \\
 \hline
  & \lambda & \dots & \lambda & 1 - K \lambda & -\frac{\xi}{2\lambda} & \frac{\xi}{2\lambda}
  \label{tab:Butcher_PFE_on_the_flyI}
\end{array}
\end{equation}

Checking the order conditions from section \ref{sec:PI_as_RK} indeed reveals second-order accuracy, independent of the outer time step size, inner time step size, and number of inner time steps. In comparison to the OPFE method, second order is achieved here using two additional function evaluations.

However, in the limit $\lambda \rightarrow 0$ using $\xi = 1$, the inner stages as well as the error estimation collapse and the IPFE scheme \eqref{tab:Butcher_PFE_on_the_fly} will degenerate to a simple Forward Euler scheme, given by the Butcher tableau
\begin{equation}
\begin{array}{ c | c }
 0 & 0 \\
 \hline
  & 1
  \label{tab:Butcher_FE}
\end{array}
\end{equation}

\section{Analysis}
\label{sec:analysis}
For analytical comparison of the schemes derived an discussed in this paper, we consider two important properties: accuracy and stability.

\subsection{Accuracy}
The accuracy of a scheme is mainly decided by the order of the scheme and its leading error coefficient. A typical first order scheme has a leading error term at time $t_{n+1}$ of the form \eqref{eq:leading_error}
\begin{equation}
    e_{n+1} = -\xi \frac{\Delta t^2}{2} u''\left(t_{n+1}\right) + \mathcal{O}\left(\Delta t^3\right),
\end{equation}

\begin{figure}[htb!]
    \centering
    \includegraphics[width=0.75\linewidth]{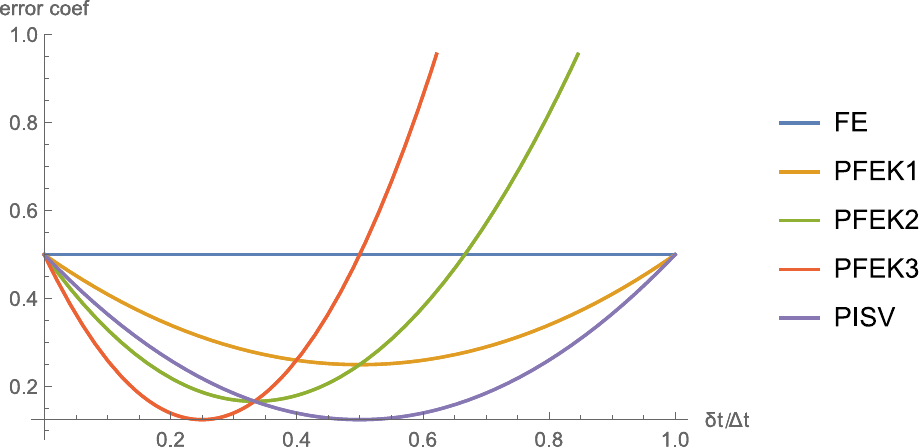}
    \caption{Error coefficient of first order schemes for FE, PFE using $K=1,2,3$ and PISV using $K=1$. Also compare table \ref{tab:leading_coefs} (left).}
    \label{fig:error_coef_first_order_methods}
\end{figure}

The second order leading error terms of Runge-Kutta methods vanish if the second order accuracy condition is fulfilled, i.e., $\sum_{j=1}^{S} b_{j}c_j = \frac{1}{2}$. This means that the mismatch $\frac{1}{2} - \sum_{j=1}^{S} b_{j}c_j$ quantifies the leading error coefficient. This yields a simple criterion to check and compare the accuracy of the previously derived schemes. Table \ref{tab:leading_coefs} shows the results for the schemes of this paper.
\begin{table}[h!]
\centering
    \begin{tabular}{ |c|c||c|c|}
         \hline
         \multicolumn{2}{|c||}{First order methods} & \multicolumn{2}{|c|}{Asymptotically second order methods} \\
         \hline
         FE   & $\frac{1}{2}$   & PRK4 ($K=1$) & $\lambda^2$\\
         PFE  & $\frac{1}{2} - K \lambda + \frac{K^2+K}{2}\lambda^2$ & PRK4 ($K=2$) & $-\frac{\lambda}{2} + 3 \lambda^2$\\
         PFE ($K=1$) & $\frac{1}{2} - \lambda + \lambda^2$  & EPHPFE ($K=2$) &  $-\frac{\lambda}{2} + 3 \lambda^2$ \\
         PISV ($K=1$) & $\frac{1}{2} - \frac{3}{2} \lambda + \frac{3}{2}\lambda^2$ & POSV ($K=2$) &  $-\lambda + 3 \lambda^2$ \\
          & & OPFE, IPFE &  $0$\\
         \hline
    \end{tabular}
    \caption{Second order leading error coefficients.}
    \label{tab:leading_coefs}
\end{table}
The first order methods in table \ref{tab:leading_coefs} clearly have a remaining term of $\frac{1}{2}$ in the limit of vanishing inner step size $\lambda \rightarrow 0$, which means that the schemes are clearly only first order. However, the PISV scheme results in a significantly reduced error coefficient in the relevant domain $\lambda \in [0,1]$, as can also be seen in \Fig \ref{fig:error_coef_first_order_methods}. In comparison with the FE scheme, all projective schemes lead to a reduction in the error coefficient for small and medium $\lambda$, indicating that they yield a more accurate solution.

\begin{figure}[htb!]
    \centering
    \includegraphics[width=0.75\linewidth]{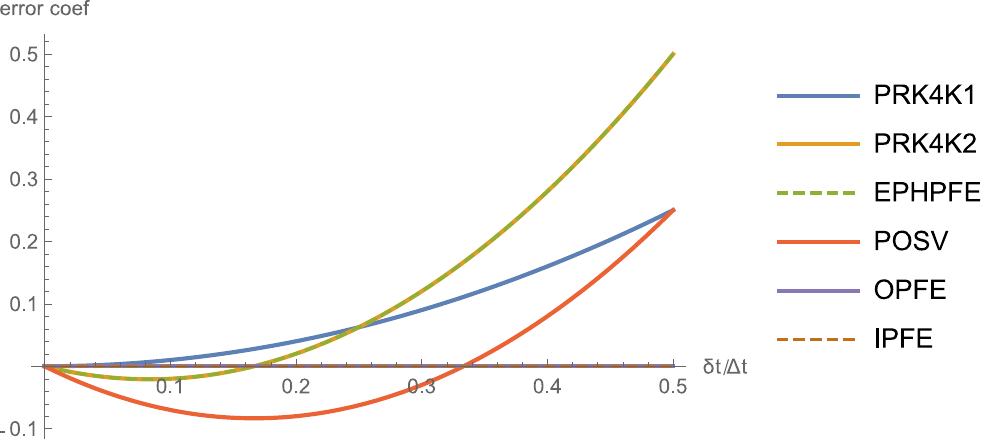}
    \caption{Error coefficient of second order schemes PRK using $K=1,2$, EPHPFE using $K=2$, POSV using $K=2$, OPFE and IPFE. Also compare table \ref{tab:leading_coefs} (right).}
    \label{fig:error_coef_second_order_methods}
\end{figure}

The methods in the right part of table \ref{tab:leading_coefs} are methods for which the second order leading erorr coefficient asymptotically vanishes for $\lambda \rightarrow 0$. This means that the scheme formally recovers second order accuracy in the case of vanishing inner step size. Notably, the PRK4 scheme and EPHPFE using $K=2$ yield the same error coefficient. In the case of OPFE and IPFE, the method is designed to yield second order regardless of the inner step size under the assumptions made in \cite{Gear2006}. The decay of the leading coefficient to zero is depicted in \Fig \ref{fig:error_coef_second_order_methods}. For very small $\lambda$, the error coefficient of the PRK4 scheme is larger when taking more inner steps $K$. This is consistent with the first order PFE scheme above.

\subsection{Stability}
\label{sec:stab_analysis}
In this section, we will see that writing projective integration schemes as RK schemes allows to use simple tools established for Runge-Kutta schemes to determine the stability properties of the schemes. Note that we are interested in stable schemes for problems including a scale separation with fast and slow modes. Therefore, we aim at stability regions that include one fast eigenvalue cluster and one (typically disconnected) slow eigenvalue cluster. We remark that schemes with connected stability regions are commonly constructed by using Telescopic Projective Integration (TPI), as covered in section \ref{sec:TPFE}. However, we focus on clearly separated fast and slow modes here.

Stability of RK schemes is commonly assessed via the scalar Dahlquist equation \cite{Hairer2000}
\begin{equation}\label{e:stab}
    u'(t) = \mu u(t), \quad u(0)= u_0,
\end{equation}
which has the exact solution $u(t) = u_0 \cdot e^{\mu t}$. For $Re(\mu) < 0$, the solution decays to zero in time. 
This decay should be mirrored by a stable numerical solution.
We apply the RK scheme to \Eqn \eqref{e:stab} and obtain $u(t + \Delta t) = g(\Delta t \mu) \cdot u(t)$, with stability function $g(\Delta t \mu) =: g(z)$.

In concise notation, the stability function of a standard RK scheme with $s$ stages is given by
\begin{equation}\label{e:stab_fcn}
    g(z) = 1 + z \vect{b}^T \left(\Vect{I} - z \Vect{A} \right)\vect{e},
\end{equation}
with $\vect{e} = (1,\ldots,1)^T \in \mathbb{R}^s$ \cite{Hairer2000}.

The behavior of the stability function determines the stability properties of the RK scheme.
A RK scheme is called stable if its solution of \eqref{e:stab} does not grow in time for $\mu < 0$, which implies $|g(z)| \leq 1$. If the scheme is stable for all $Re(z) < 0$, then it is called A-stable, which can only be fulfilled by implicit RK schemes. Explicit schemes are typically only stable in a small part of the negative half plane, i.e., for small time steps or slow modes around $z=-1$ \cite{Hairer2000}. Projective Integration schemes aim to extend the stability region such that the scheme is also stable for fast modes around $z = -\frac{1}{\epsilon}$, with $\epsilon \ll 1$ \cite{Melis2017}. In the following figures, we mark both points $z = -1$ and $z = -\frac{1}{\epsilon}$ with a black dot to indicate the desired stability region. Note that the appearance of a domain around $z = -\frac{1}{\epsilon}$ inside the stability region is a typical property of PI schemes, due to the use of an inner time stepping scheme with small time step size $\delta t = \epsilon$. More details on the stability properties of PI schemes can be found in \cite{Lafitte2016,Lafitte2017,Melis2017}.

\begin{figure}[htb!]
    \centering
    \begin{subfigures}
    \subfloat[FE.
    ]{\includegraphics[width=0.5\linewidth]{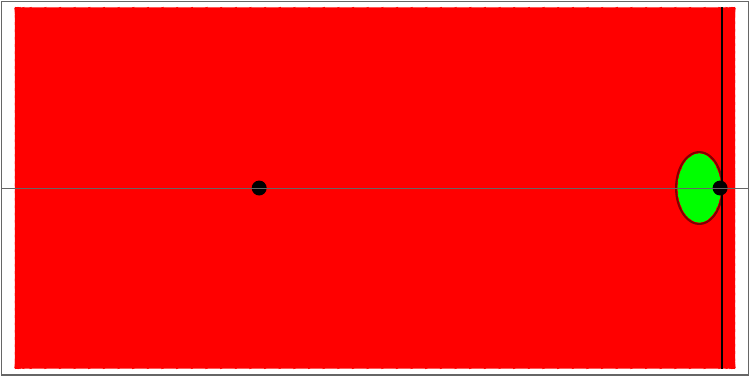}}~
    \subfloat[RK4.
    ]{\includegraphics[width=0.5\linewidth]{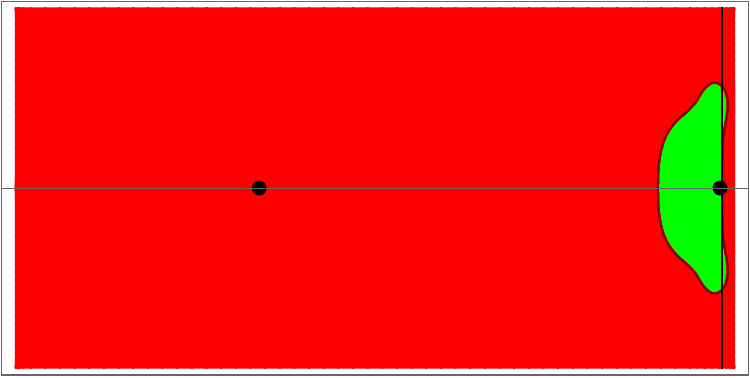}}\\
    \subfloat[PFE, $K+1=1$.
    ]{\includegraphics[width=0.5\linewidth]{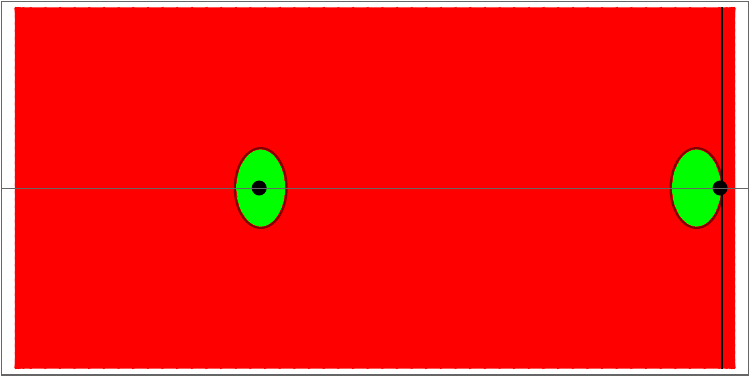}}~
    \subfloat[PRK4, $K+1=1$.
    ]{\includegraphics[width=0.5\linewidth]{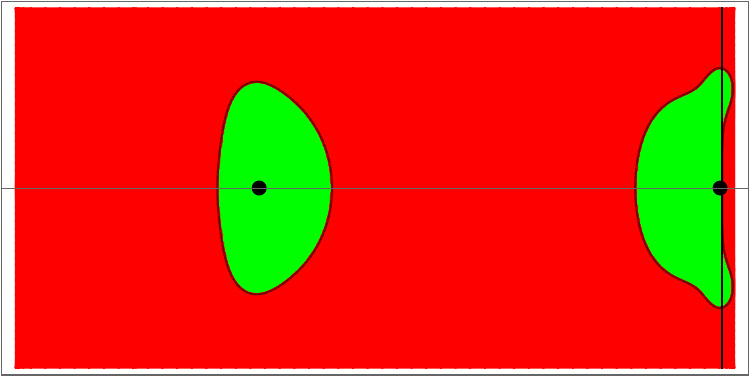}}\\
    \subfloat[PFE, $K+1=2$.
    ]{\includegraphics[width=0.5\linewidth]{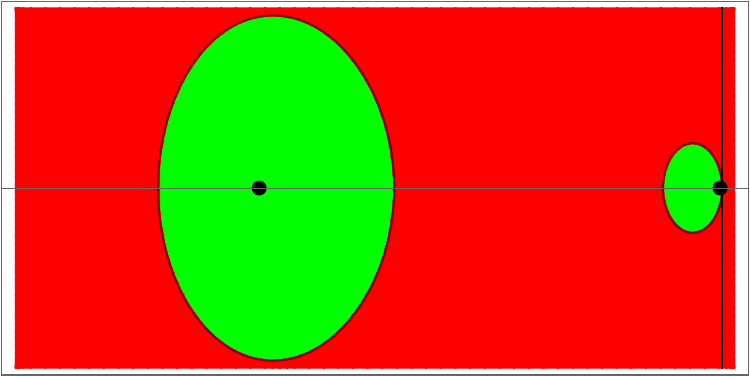}}~
    \subfloat[PRK4, $K+1=2$.
    ]{\includegraphics[width=0.5\linewidth]{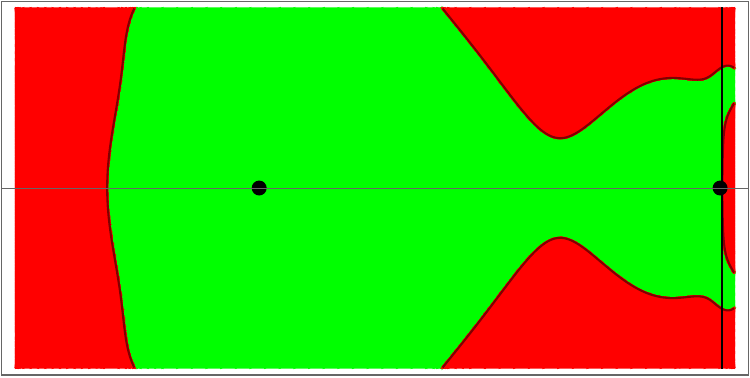}}
    \end{subfigures}
    \caption{Stability regions of FE scheme (a), RK4 scheme (b), PFE scheme using $K+1=1$ (c), PRK4 scheme using $K+1=1$ (d), PFE scheme using $K+1=2$ (e) and PRK4 scheme using $K+1=2$ (f). Black dots at $(-1,0)$ and $(-\frac{1}{\epsilon},0)$ indicate eigenvalue clusters for desired stability region. Projective schemes are also stable for fast cluster and larger $K$ increases the stability region.}
    \label{fig:stab_FE_RK4}
\end{figure}

In \Fig \ref{fig:stab_FE_RK4}, the stability regions for Forward Euler methods (FE) (left column) and Runge-Kutta 4 methods (RK4) (right column) are shown. The standard methods FE (a) and RK4 (b) are only stable for slow modes as indicated by the stability region around the slow cluster $z=-1$ near the origin. The projective variants PFE (c) and PRK4 (d) using $K+1=1$ inner iteration, however, clearly show that the stability region is augmented by a stable region around the fast cluster, while the PRK4 has a slightly increased stability region. Using $K+1=2$ inner iterations in (e) and (f), the stability regions increase and even yield a connected stability region in case of the PRK4 scheme (f).

\begin{figure}[htb!]
    \centering
    \begin{subfigures}
    \subfloat[EMR.
    ]{\includegraphics[width=0.5\linewidth]{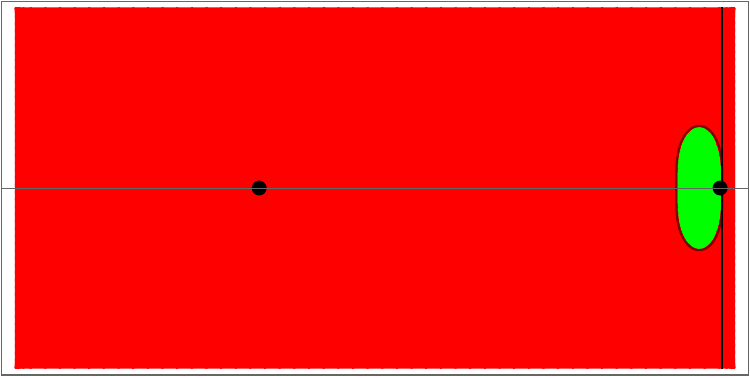}}~
    \subfloat[EPHPFE, $K+1=2$.
    ]{\includegraphics[width=0.5\linewidth]{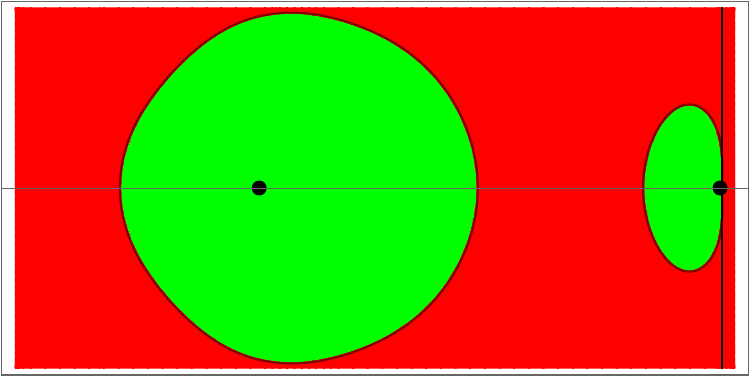}}\\
    \subfloat[PISV, $K+1=1$.
    ]{\includegraphics[width=0.5\linewidth]{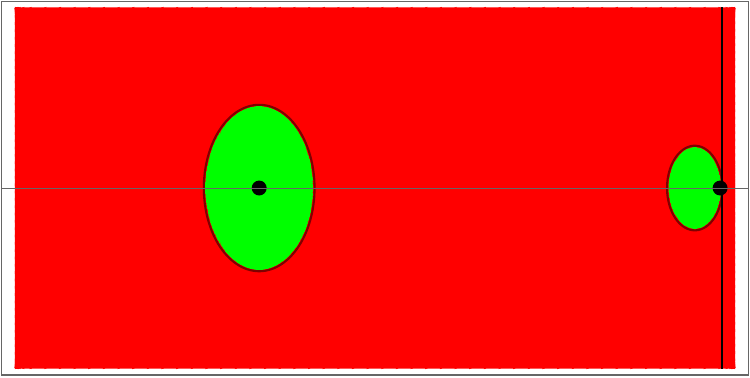}}~
    \subfloat[POSV, $K+1=2$.
    ]{\includegraphics[width=0.5\linewidth]{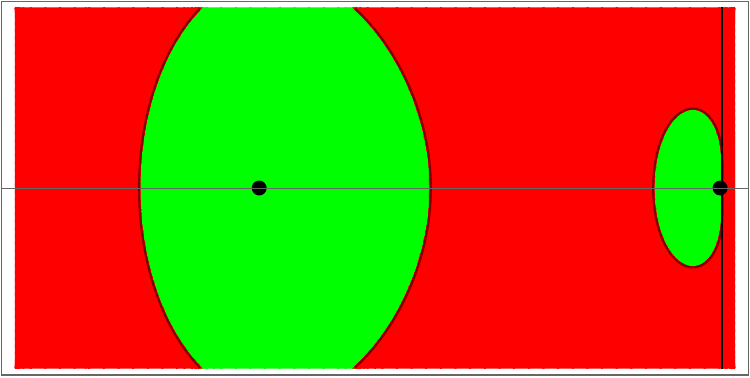}}\\
    \subfloat[OPFE, $K+1=2$.
    ]{\includegraphics[width=0.5\linewidth]{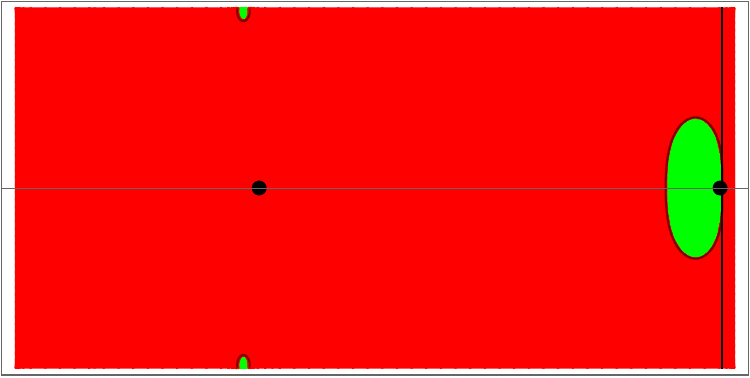}}~
    \subfloat[IPFE, $K+1=3$.
    ]{\includegraphics[width=0.5\linewidth]{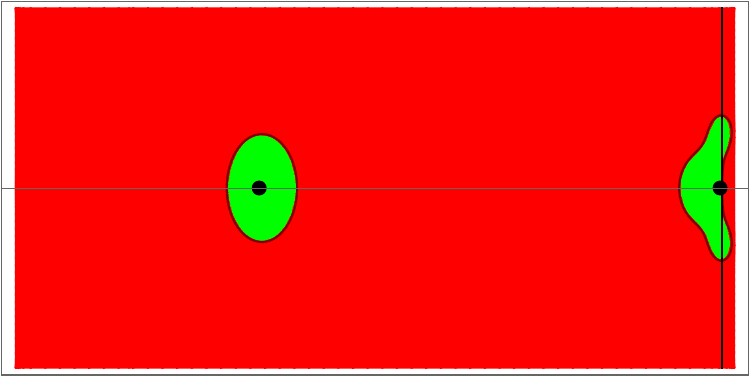}}
    \end{subfigures}
    \caption{Stability region of EMR scheme (a), EPHPFE scheme using $K+1=2$ (b), PISV scheme using $K+1=1$ (c), POSV scheme using $K+1=2$ (d), OPFE scheme using $K+1=2$, and IPFE scheme using $K=3$. Black dots at $(-1,0)$ and $(-\frac{1}{\epsilon},0)$ indicate eigenvalue clusters for desired stability region. Projective embedded schemes and IPFE are stable for fast cluster and show larger stability region for larger $K$.}
    \label{fig:stab_new_schemes}
\end{figure}

In \Fig \ref{fig:stab_new_schemes}, the stability regions for the new time adaptive embedded schemes derived in the previous sections are plotted. The standard explicit midpoint rule (EMR) in (a), in contrast, only results in stable integration of the slow cluster whereas the Embedded Projective Heun Projective Forward Euler scheme (EPHPFE) from \Eqn \ref{tab:Butcher_ePRKK2_cor} in (b) yields stable integration of the fast cluster, too. The same holds true for the scheme with projective outer step size variation (POSV) from Section \ref{sec:POSV} in (d) as well as for the scheme with projective inner step size variation (PISV) from Section \ref{sec:PISV} in (c). Notably, the stability region of the POSV scheme seems to be larger than that of the PISV scheme. The POSV scheme results in a stability region close to the PFE scheme with $K+1=2$ from \Fig \ref{fig:stab_FE_RK4}.

The two method using on-the-fly error estimation in \Fig \ref{fig:stab_new_schemes}(e) and (f) require some extra explanation. It can be seen that the OPFE scheme does not yield a stable integration of the fast cluster. This is due to the derivative estimation using the outer integration points, which are unable to capture the fast dynamics of the inner iterations. The IPFE method in \Fig \ref{fig:stab_new_schemes} (f) removes that disadvantage by using an additional inner iteration to estimate the derivative and therefore captures the fast dynamics. This leads to a stable integration of the fast eigenvalue cluster. In other tests (not shown for conciseness) we could show that the stability region of the IPFE increases with increasing $K$, especially around the fast eigenvalue cluster.

Summarizing, we see that the new projective schemes indeed allow for stable integration of fast clusters due to the extended stability region of the schemes. This can be analyzed straightforwardly in the framework of RK schemes.

\section{Numerical tests}
\label{sec:conv_test}
\begin{figure}[htb!]
    \centering
    \begin{subfigures}
    \subfloat[$\epsilon=0.05$.
    ]{\includegraphics[width=0.45\linewidth]{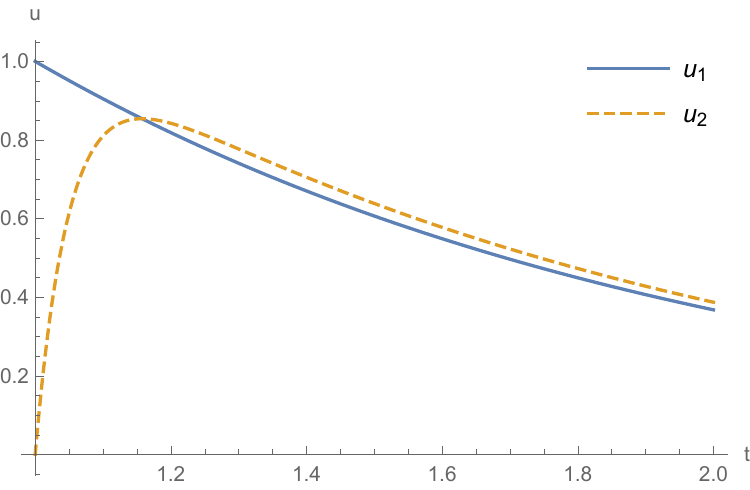}}~~
    \subfloat[$\epsilon=0.005$.
    ]{\includegraphics[width=0.45\linewidth]{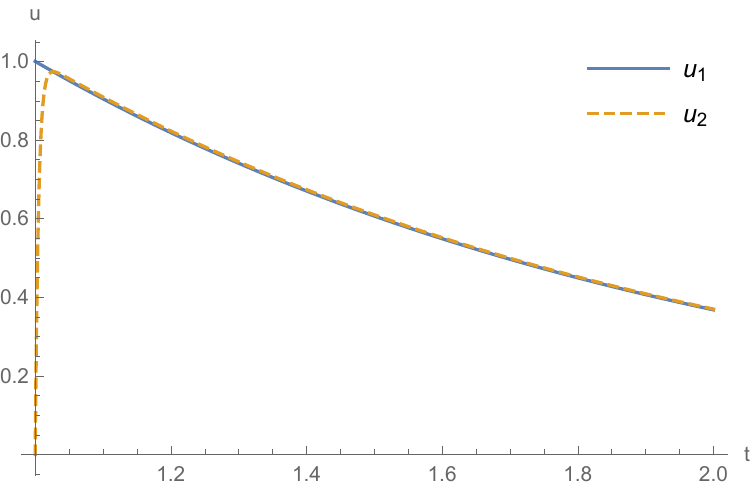}}
    \end{subfigures}
    \caption{Exact solution of test case \eqref{e:conv_test}-\eqref{e:conv_test2} with varying $\epsilon=0.05$ (left) and $\epsilon=0.005$ (right). Note that $u_2$ relaxes towards $u_1$ much faster with decreasing $\epsilon$, leading to a stiff model.}
    \label{fig:testcase}
\end{figure}

While the focus of this work is the derivation of Projective Integration schemes in the framework of RK schemes, it is important to access the numerical properties of the new schemes. We therefore perform a simple convergence test to investigate the speed of convergence for the following two-scale model problem consisting of two coupled equations with separate scales
\begin{eqnarray}
    u_1'(t) &=& - \alpha u_1(t),                     \quad\quad\quad\quad\quad u_1(0)= 1, \label{e:conv_test} \\
    u_2'(t) &=& \frac{1}{\epsilon} \left(u_1(t) - u_2(t) \right),   \quad u_2(0)= 0,\label{e:conv_test2}
\end{eqnarray}
where we use $\alpha = 1$ and $\epsilon \ll 1$. Note that the Jacobian of the model problem \eqref{e:conv_test}-\eqref{e:conv_test2} reads
\begin{equation}
    A = \left(
          \begin{array}{cc}
            -\alpha & 0 \\
            \frac{1}{\epsilon} & -\frac{1}{\epsilon} \\
          \end{array}
        \right),
\end{equation}
with eigenvalues $\lambda_1 = -\alpha$, $\lambda_2 = -\frac{1}{\epsilon}$, leading to $u_1$ and $u_2$ developing on different scales. The exact solution is given by
\begin{eqnarray}
    u_1(t) &=& u_1(0) \exp(-\alpha t),    \\
    u_2(t) &=& \frac{u_1(0)}{1- \alpha\epsilon} \exp(-\alpha t)  + \frac{(1-\alpha \epsilon) u_2(0) - u_1(0)}{1- \alpha\epsilon} \exp\left(-\frac{t}{\epsilon}\right) ,
\end{eqnarray}
so that $u_1$ slowly relaxes to zero in time while $u_2$ approaches $u_1$ with potentially fast relaxation time $\frac{1}{\epsilon}$. The exact solution is plotted for $\epsilon = 0.05$ and $\epsilon = 0.005$ in Figure \ref{fig:testcase}. This test case is suitable for assessing the accuracy and stability properties of the time integration methods proposed above as it includes a scale separation and requires efficient techniques to speed up standard methods.

Standard FE or RK4 schemes will not be stable for the integration of the model problem \eqref{e:conv_test}-\eqref{e:conv_test2} for time steps $\Delta t > \epsilon$, due to the fast evolution of $u_2$. We will subsequently assess the numerical convergence of the new projective schemes and the error estimators of the embedded projective schemes. For a concise presentation of the results, we show error plots for $u_1$. However, for large time $t$, the solution will yield $u_1(t) \approx u_2(t)$ due to the relaxation in \eqref{e:conv_test2}, so that the errors in our numerical simulations where practically the same for $u_1$ and $u_2$.

\begin{remark}(\textit{Other test cases})
    While this paper concerns the derivation and analysis of PI schemes in the framework of RK methods and a simple but useful test case is employed, other models from applications in science and engineering can readily be applied to the methods derived in this paper.
    This includes models from non-equilibrium rarefied gas dynamics \cite{Torrilhon2016}, where the macroscopic variables density $\rho$, bulk velocity $v$ and temperature $T$ evolve on a macroscopic fluid scale, while non-equilibrium variables potentially evolve much faster. The spectrum can contain clear spectral gaps as investigated in \cite{Koellermeier2021,Koellermeier2020h}.
    Another application can be models for extended shallow water equations \cite{Kowalski2019}, which are potentially stiff with multiple scales as investigated in \cite{Amrita2022,Huang2022}. Investigation of spatially adaptive methods based on \cite{KoellermeierAPI} might be beneficial for additional computational speedup.
\end{remark}

\subsection{Numerical error convergence}
As first of two numerical tests, we investigate the numerical error convergence of the method by monitoring the numerical error with respect to the analytical solution for vanishing $\Delta t \rightarrow 0$ and constant value of $\lambda = \frac{\delta t}{\Delta t}$, where as usual the inner time step size is taken as $\delta t = \epsilon$.

\begin{figure}[htb!]
    \centering
    \begin{subfigures}
    \subfloat[Decreasing $\Delta t \rightarrow 0$, with constant $\delta t / \Delta t$.
    ]{\includegraphics[width=0.7\linewidth]{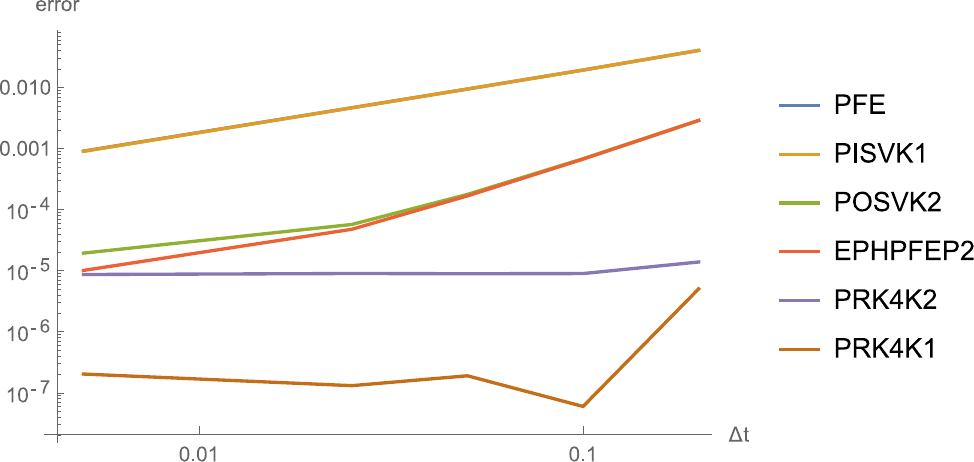}}\\
    \subfloat[Decreasing $\Delta t \rightarrow 0$ and $\delta t / \Delta t \rightarrow 0$.
    ]{\includegraphics[width=0.7\linewidth]{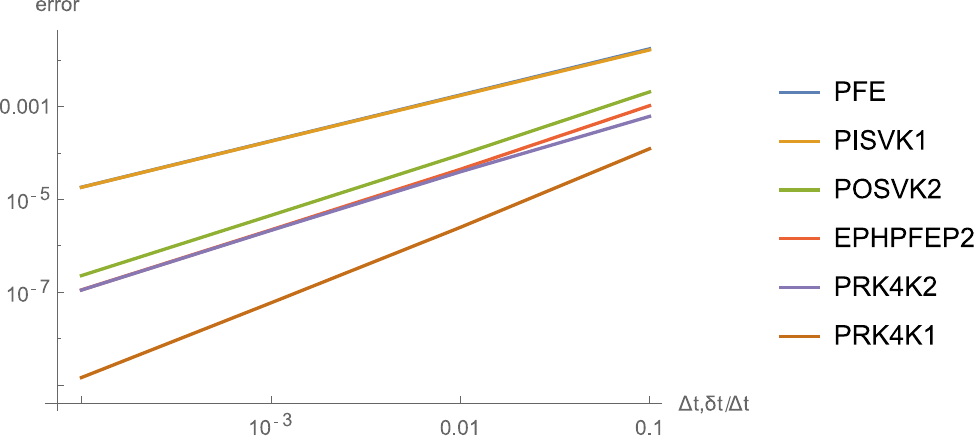}}
    \end{subfigures}
    \caption{Numerical error convergence of schemes for decreasing time step size $\Delta t$. If only $\Delta t$ goes to zero (a), convergence is not obtained because at some point the $\lambda = \delta t / \Delta t$ error becomes dominant (a). If also $\lambda = \delta t / \Delta t$ goes to zero (b), the expected convergence is obtained.}
    \label{fig:stab_conv}
\end{figure}
We choose $\epsilon = 10^{-5}$ and see in \Fig \ref{fig:stab_conv}(a) that the PFE ($K+1=1$), PISV ($K+1=1$), POSV ($K+1=2$), and EPHPFE ($K+1=2$) schemes seem to correctly converge with first order. Both PRK4 schemes for $K+1=1$ and $K+1=2$, however, do not seem to converge for smaller values of $\Delta t$. This is due to the increasing value of $\lambda = \frac{\delta t}{\Delta t}$, which can no longer be neglected for small $\Delta t$. We therefore design a second test were both $\Delta t \rightarrow 0$ and $\lambda = \frac{\delta t}{\Delta t} \rightarrow 0$. This is done by decreasing the value of $\epsilon = \delta t$ simultaneously. In this situation, we see in \Fig \ref{fig:stab_conv}(b) that all errors decrease and the schemes obtain the expected convergence. Note that the schemes are formally only of first order as shown in Theorem \ref{th:consistency_PRK}. However, we clearly see that the PR4 schemes are the most accurate, while the simple PFE and PISV schemes yield the largest errors.

The numerical convergence test therefore shows that the projective integration schemes derived in this paper obtain the proven convergence rates from Theorem \ref{th:consistency_PRK} while more complex PRK schemes are still beneficial to obtain more accurate solutions.

\subsection{Embedded schemes and error estimators}
A numerical study was conducted to assess the performance of the error estimators and the resulting error behavior of the newly derived embedded projective schemes. We consider the test problem \eqref{e:conv_test}-\eqref{e:conv_test2}, for stable settings $\Delta t = 0.1$, $\delta t = \epsilon$ and vary $\delta t / \Delta t \in (0,1]$ to investigate the error behavior depending on the inner time step size. We only consider the error of the first equation \eqref{e:conv_test} modeling exponential decay after one standard time step $\Delta t$ and are interested in the following four quantities per embedded scheme:
\begin{itemize}
  \item The error of the corrected embedded scheme: $\parallel u_{\ast}^{n+1} - u(t_{n+1}) \parallel $ (in blue, if available)
  \item The error of the lower-order solution error: $\parallel \widetilde{u}^{n+1} - u(t_{n+1}) \parallel$ (in yellow)
  \item The error of the higher-order solution: $\parallel u^{n+1} - u(t_{n+1}) \parallel$ (in green)
  \item The error estimate of the embedded method: $\parallel u^{n+1} - \widetilde{u}^{n+1} \parallel$ (in red)
\end{itemize}
The numerical errors of the new embedded projective schemes will be compared with the respective standard non-projective embedded schemes (if available).

\begin{figure}[htb!]
    \centering
    \begin{subfigures}
    \subfloat[EMR.
    ]{\includegraphics[width=0.5\linewidth]{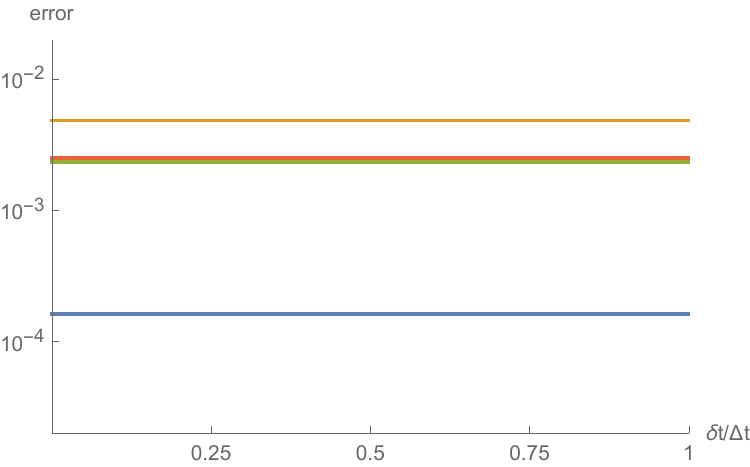}}~
    \subfloat[POSV, $K+1=1$.
    ]{\includegraphics[width=0.5\linewidth]{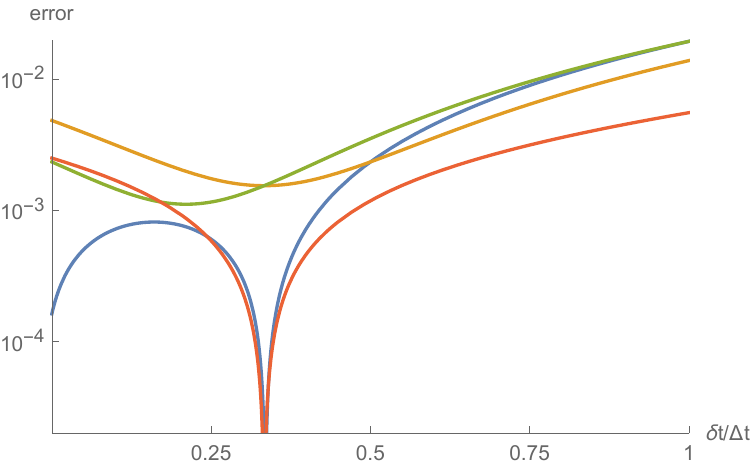}}\\
    \subfloat[EHFE.
    ]{\includegraphics[width=0.5\linewidth]{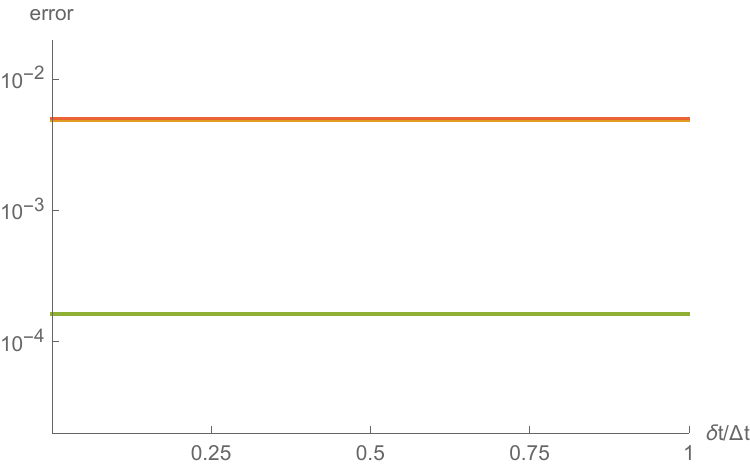}}~
    \subfloat[EPHPFE, $K+1=2$.
    ]{\includegraphics[width=0.5\linewidth]{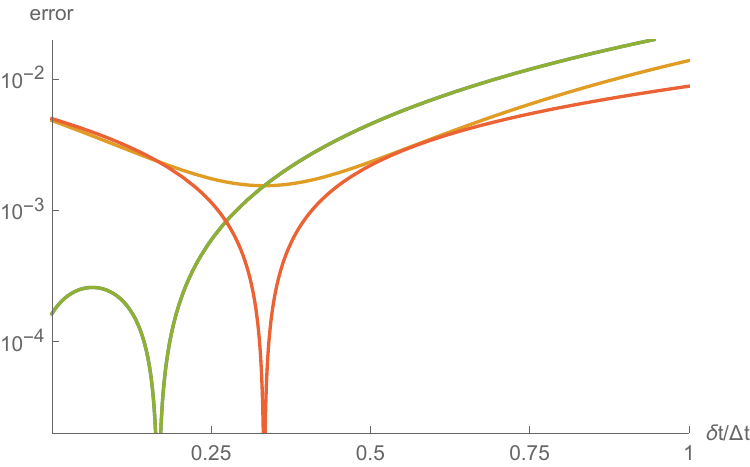}}\\
    \subfloat[PISV, $K+1=2$.
    ]{\includegraphics[width=0.8\linewidth]{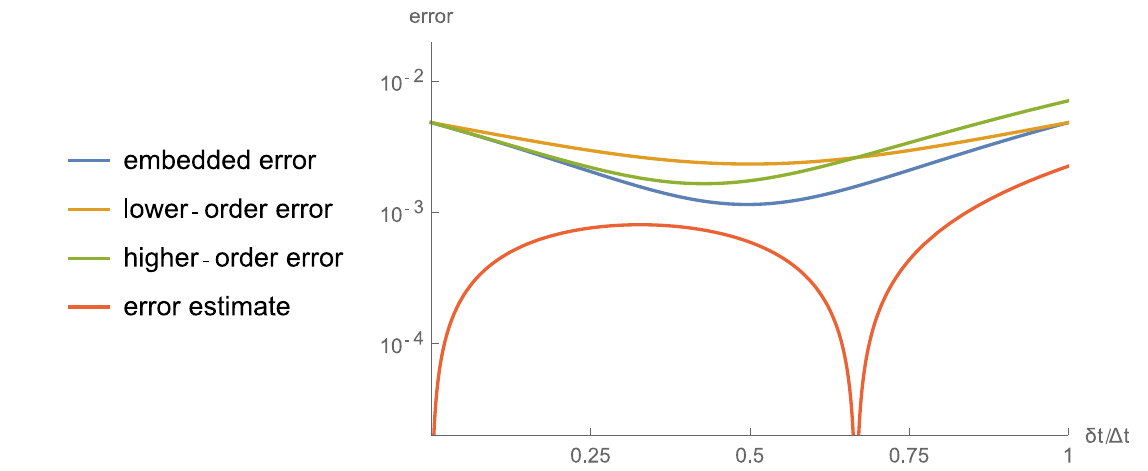}}
    \end{subfigures}
    \caption{Errors and error estimators of embedded projective schemes (b), (d), (e) in comparison to standard embedded schemes (a), (c) for changing inner time step size $\delta t$. 
    Full embedded solution error (blue), lower-order solution error (yellow), higher-order solution (green), error estimate (red). The projective versions (b),(d) show smaller errors than the standard schemes (a), (c), resulting in smaller error estimates.}
    \label{fig:errors}
\end{figure}
In Figure \ref{fig:errors} we see the respective data for the embedded schemes also investigated in Figure \ref{fig:stab_new_schemes}. On the one hand, the standard Explicit Midpoint Rule (EMR) in \ref{fig:errors}(a) does not contain inner time steps and therefore does not depend on the inner time step size $\delta t$. Note, how the error estimator (in red) measures the error of the higher-order solution (in green). On the other hand, the new POSV scheme in \ref{fig:errors}(b) consistently shows the same error values as the EMR scheme for vanishing $\delta t / \Delta t \rightarrow 0$, but it has a decreased error for intermediate values of $\delta t / \Delta t$. For values $\delta t / \Delta t = \mathcal{O}(1)$, the scheme is less accurate again, but a projective scheme would not be used in that case.

Similarly, the embedded Heun Forward Euler (EHFE) method in Figure \ref{fig:errors}(c) does not contain inner time steps and therefore does not depend on the inner time step size $\delta t$, too. Note, how the error estimator (in red) here measures the error of the low-order solution (in yellow). Correspondingly, the new projective version EPHPFE in \ref{fig:errors}(d) shows the same error values as the EHFE scheme for vanishing $\delta t / \Delta t \rightarrow 0$. However, for non-vanishing $\delta t / \Delta t$ the new scheme accurately predicts the error using its build-in error estimator. For unusually large $\delta t / \Delta t $, the error is slightly increased again, while the error estimator shows good performance in measuring the lower-order error.

Lastly, we see that also the new PISV method (for which there is no corresponding standard scheme, since it is based on varying the inner step size) shows a comparable behavior with decreasing error for intermediate values of $\delta t / \Delta t$. The error of the corrected embedded scheme (in blue) is clearly smaller than both the lower-order and higher-order solutions for all values $\delta t \in (0,\Delta t)$. However, since the method is based on inner time step variation, the error estimate does not give an accurate assessment of the actual error. This is expected as the inner time step size is typically not chosen because of accuracy, but because of stability constraints. The PISV method is less advantageous than the POSV method with respect to accurate error control.

Summarizing, the new embedded projective schemes POSV, EPHPFH, and POSV show a favorable error behavior in comparison to their corresponding standard schemes for small and intermediate inner time step sizes, which is exactly the application case for projective schemes. This indicates a good performance of the respective error estimators of these new embedded projective schemes.

\section{Conclusion}
In this paper, we use the definition of Projective Integration (PI) methods as explicit time stepping schemes to write them as Runge-Kutta (RK) methods. This remarkably simple rewriting allows to check their consistency and accuracy properties easily by means of the order conditions of the RK method, without performing tedious Taylor expansions, which is especially troublesome for the many steps of a Projective Runge-Kutta (PRK) method or Telescopic Projective Integration (TPI) method. Spatially adaptive Projective Integration methods can be included as partitioned RK methods. Using the framework of RK methods, we derived new time adaptive PI schemes, that are based on the corresponding embedded RK method, step size variation, or an on-the-fly error estimation.

The accuracy, stability, and numerical convergence properties of the rewritten methods and newly derived methods were easily analyzed in the RK framework and the projective methods clearly show the desired enlarged stability region while converging with the theoretically derived order of accuracy. We prove instability of an on-the-fly error estimation method and show stability of a new improved version using the rewriting as RK method.

The work in this paper allows to use more tools for stability analysis of RK schemes to be applied to PI.
The number and size of inner time steps of PI methods or the different parameters of a TPI method could then be optimized to answer the question if similarly stable RK schemes with fewer stages can be derived. 
The development of space \emph{and} time adaptive PI methods is another possible future research direction.\\

\section*{Acknowledgements}
This research has been partially supported by the European Union’s Horizon 2020 research and innovation program under the Marie Sklodowska-Curie grant agreement no. 888596.
The authors would like to acknowledge the financial support of the CogniGron research center and the Ubbo Emmius Funds (University of Groningen).

\section*{Conflict of interest statement}
The authors certify that they have no affiliations with or involvement in any organization or entity with any conflicting interest in the subject matter or materials discussed in this manuscript.

\section*{Data availability statement}
The datasets generated during and/or analysed during the current study are available from the corresponding author on reasonable request

\appendix

\section{Examples: PRK4 schemes}
\label{app:PRK4}
The Butcher tableaus for the projective versions of the standard RK scheme of fourth order using different inner step numbers are given below.

For the PRK4K1 version with $K=1$, i.e., $K+1=2$ inner steps, the Butcher tableau reads:
\begin{equation}
\begin{array}{ l | c c c c c c c c}
 0                      & 0       &                         &           &               &           &               &         &     \\
 0 + \lambda            & \lambda & 0                       &           &               &           &               &         &     \\
 \frac{1}{3}            & \lambda & \frac{1}{3} - \lambda   & 0         &               &           &               &         &     \\
 \frac{1}{3} + \lambda  & \lambda & \frac{1}{3} - \lambda   & \lambda   & 0             &           &               &         &     \\
 \frac{2}{3}            & \lambda & -\frac{1}{3} + 2\lambda & 0         & 1 - 3\lambda  & 0         &               &         &     \\
 \frac{2}{3} + \lambda  & \lambda & -\frac{1}{3} + 2\lambda & 0         & 1 - 3\lambda  & \lambda   & 0             &         &     \\
 1                      & \lambda & 1 - \lambda             & 0         & -1 + 2\lambda & 0         & 1 - 2 \lambda & 0       &     \\
 1 + \lambda            & \lambda & 1 - \lambda             & 0         & -1 + 2\lambda & 0         & 1 - 2 \lambda & \lambda & 0   \\
 \hline
                        & \lambda & \frac{1}{8} + \frac{3}{4}\lambda & 0 & \frac{3}{8} - \frac{6}{8}\lambda & 0 & \frac{3}{8} - \frac{6}{8}\lambda & 0 & \frac{1}{8} - \frac{1}{4}\lambda
 \label{tab:Butcher_PRK4K1}
\end{array}
\end{equation}

The PRK4K2 version with $K=2$, i.e., $K+1=3$ inner steps, inner steps is given by the following Butcher tableau:
\begin{equation}
\begin{array}{ l | c c c c c c c c c c c c}
 0                      & 0       &         &                                 &           &               &                       &           &               &               &         &           \\
 0 + \lambda            & \lambda & 0       &                                 &           &               &                       &           &               &               &         &           \\
 0 + 2\lambda           & \lambda & \lambda & 0                               &           &               &                       &           &               &               &         &           \\
 \frac{1}{3}            & \lambda & \lambda & \frac{1}{3} -2\lambda           & 0         &               &                       &           &               &               &         &           \\
 \frac{1}{3} + \lambda  & \lambda & \lambda & \frac{1}{3} -2\lambda           & \lambda   & 0             &                       &           &               &               &         &           \\
 \frac{1}{3} +2\lambda  & \lambda & \lambda & \frac{1}{3} -2\lambda           & \lambda   &  \lambda      & 0                     &           &               &               &         &           \\
 \frac{2}{3}            & \lambda & \lambda & -\frac{1}{3}-\frac{5}{2}\lambda & 0         & 0             & 1-\frac{9}{2}\lambda  & 0         &               &               &         &           \\
 \frac{2}{3} + \lambda  & \lambda & \lambda & -\frac{1}{3}-\frac{5}{2}\lambda & 0         & 0             & 1-\frac{9}{2}\lambda  & \lambda   & 0             &               &         &           \\
 \frac{2}{3} +2\lambda  & \lambda & \lambda & -\frac{1}{3}-\frac{5}{2}\lambda & 0         & 0             & 1-\frac{9}{2}\lambda  & \lambda   & \lambda       & 0             &         &           \\
 1                      & \lambda & \lambda & 1 - 2\lambda                    & 0         & 0             & -1 + 3\lambda         & 0         & 0             & 1 - 3 \lambda & 0       &           \\
 1 + \lambda            & \lambda & \lambda & 1 - 2\lambda                    & 0         & 0             & -1 + 3\lambda         & 0         & 0             & 1 - 3 \lambda & \lambda & 0         \\
 1 +2\lambda            & \lambda & \lambda & 1 - 2\lambda                    & 0         & 0             & -1 + 3\lambda         & 0         & 0             & 1 - 3 \lambda & \lambda & \lambda   \\
 \hline
                        & \lambda & \lambda & \frac{1}{8} + \frac{5}{8}\lambda & 0 & 0 & \frac{3}{8} - \frac{9}{8}\lambda & 0 & 0 & \frac{3}{8} - \frac{9}{8}\lambda & 0 & 0 & \frac{1}{8} - \frac{3}{8}\lambda
 \label{tab:Butcher_PRK4K2}
\end{array}
\end{equation}

\bibliographystyle{plain}
\bibliography{library_fixed}

\end{document}